%% file: main.tex
\documentclass[a4paper, 12pt]{article}
\pdfoutput=1

\input{_preamble}

\input{_commands}

\title{A motivic integral identity\\for \titleminusone-shifted symplectic stacks}
\author{Chenjing Bu}
\date{}

\begin{document}

\initlengths

\maketitle

\begin{abstract}
    \input{abstract}
\end{abstract}

{
    \hypersetup{linkcolor=black}
    \setstretch{1.2}
    \tableofcontents
}

\clearpage
\section{Introduction}

\input{intro}

\clearpage
\section{Motivic vanishing cycles}

\input{mot}

\clearpage
\section{Graded and filtered points}

\input{grad}

\clearpage
\section{The main results}

\input{id}

\newpage
\phantomsection
\addcontentsline{toc}{section}{References}
\sloppy
\setstretch{1.1}
\renewcommand*{\bibfont}{\normalfont\small}
\setlength{\bibitemsep}{0.2\baselineskip}
\printbibliography

\par\noindent\rule{0.38\textwidth}{0.4pt}
{\par\noindent\small
\hspace*{2em}Chenjing Bu\hspace{4em}\texttt{bu@maths.ox.ac.uk}\\[-2pt]
\hspace*{2em}Mathematical Institute, University of Oxford, Oxford OX2 6GG, United Kingdom.
}

\end{document}

%% file: _preamble.tex
\usepackage[margin=25mm]{geometry}

\usepackage[utf8]{inputenc}
\usepackage[T1]{fontenc}
\usepackage[british]{babel}
\usepackage[style=british]{csquotes}
\usepackage[en-GB]{datetime2}
\DTMlangsetup[en-GB]{ord=omit}
\usepackage{xcolor}

\usepackage{amsthm, thmtools}
\usepackage{mathtools}
\usepackage{titlesec} 
\usepackage{hyperref}
\hypersetup{
    bookmarksnumbered,
    colorlinks,
    linkcolor={cyan!50!blue},
    citecolor={cyan!50!blue},
    urlcolor={cyan!50!blue}
}
\usepackage[capitalise,nameinlink]{cleveref}
\usepackage[shortcuts]{extdash}

\usepackage{tikz}
\usepackage{tikz-cd}

\usepackage[lining,tabular,scale=.9]{FiraSans}
\usepackage[tt=false,semibold]{libertine}
\usepackage[libertine,timesmathacc,smallerops]{newtxmath}
\usepackage[lining,scaled=.8]{FiraMono}
\usepackage[cal=boondox,frak=euler]{mathalpha}

\DeclareFontEncoding{MDA}{}{}
\DeclareSymbolFont{mathdesignA}{MDA}{mdput}{m}{n}
\SetSymbolFont{mathdesignA}{bold}{MDA}{mdput}{b}{n}
\DeclareSymbolFontAlphabet{\mathbb}{mathdesignA}

\tikzcdset{
    arrow style=tikz,
    arrows={/tikz/line width=.5pt},
    diagrams={>={Straight Barb[scale=0.8]}}
}

\DeclareFontFamily{OMX}{MnSymbolE}{}
\DeclareSymbolFont{MnLargeSymbols}{OMX}{MnSymbolE}{m}{n}
\SetSymbolFont{MnLargeSymbols}{bold}{OMX}{MnSymbolE}{b}{n}
\DeclareFontShape{OMX}{MnSymbolE}{m}{n}{
    <-6>  MnSymbolE5
   <6-7>  MnSymbolE6
   <7-8>  MnSymbolE7
   <8-9>  MnSymbolE8
   <9-10> MnSymbolE9
  <10-12> MnSymbolE10
  <12->   MnSymbolE12
}{}
\DeclareFontShape{OMX}{MnSymbolE}{b}{n}{
    <-6>  MnSymbolE-Bold5
   <6-7>  MnSymbolE-Bold6
   <7-8>  MnSymbolE-Bold7
   <8-9>  MnSymbolE-Bold8
   <9-10> MnSymbolE-Bold9
  <10-12> MnSymbolE-Bold10
  <12->   MnSymbolE-Bold12
}{}
\DeclareMathDelimiter{[}{\mathopen}{MnLargeSymbols}{'000}{MnLargeSymbols}{'000}
\DeclareMathDelimiter{]}{\mathclose}{MnLargeSymbols}{'005}{MnLargeSymbols}{'005}
\DeclareMathDelimiter{\llbr}{\mathopen}{MnLargeSymbols}{'102}{MnLargeSymbols}{'102}
\DeclareMathDelimiter{\rrbr}{\mathclose}{MnLargeSymbols}{'107}{MnLargeSymbols}{'107}

\usepackage{bm}
\usepackage{cases}
\usepackage{accents}

\usepackage{setspace}
\usepackage{subdepth} 

\newcommand{\initlengths}{%
    \setlength{\abovedisplayshortskip}{3pt plus 9pt minus 3pt}%
    \setlength{\belowdisplayshortskip}{9pt plus 9pt minus 9pt}%
    \setlength{\abovedisplayskip}{9pt plus 9pt minus 9pt}%
    \setlength{\belowdisplayskip}{9pt plus 9pt minus 9pt}%
    \tolerance 400
}

\usepackage{titling}
\pretitle{\vspace{-36pt}\setstretch{1.1}\begin{center}\LARGE\libertineDisplay\fontdimen2\font=0.3em} 
\posttitle{\end{center}\vspace{3pt}}

\numberwithin{paragraph}{subsection}
\newcommand{\parafont}{\bfseries}
\newcommand{\parasep}{9pt plus 3pt minus 3pt}

\titleformat{\section}{\Large\libertineDisplay}{\thesection}{1em}{}
\titleformat{\subsection}{\large\firamedium\boldmath}{\thesubsection}{1em}{}
\titleformat{\paragraph}[runin]{\parafont}{\theparagraph.}{.33em}{\normalfont\bfseries\boldmath}
\titlespacing{\paragraph}{0pt}{\parasep}{.5em}

\usepackage{tocloft}

\renewenvironment{abstract}{%
    \centering\begin{minipage}{.85\textwidth}%
    \setlength{\parindent}{1.5em}%
    \centerline{\firamedium\abstractname}%
    \setstretch{1.15}%
    \par\vspace{6pt}%
}{\end{minipage}\par\vspace{18pt}}

\makeatletter
\declaretheoremstyle[
    spaceabove=\parasep, spacebelow=\parasep,
    postheadspace=.5em,
    headfont=\normalfont\bfseries,
    headpunct={},
    headformat={\NUMBER.\@\ \NAME.\@\NOTE},
    notefont=\normalfont\bfseries\boldmath,
    notebraces={}{.},
    bodyfont=\itshape,
]{theorem}
\declaretheoremstyle[
    spaceabove=\parasep, spacebelow=\parasep,
    postheadspace=.5em,
    headfont=\normalfont\bfseries,
    headpunct={},
    headformat={\NAME.\@\NOTE},
    notefont=\normalfont\bfseries\boldmath,
    notebraces={}{.},
    bodyfont=\itshape,
]{theorem*}
\declaretheoremstyle[
    spaceabove=\parasep, spacebelow=\parasep,
    postheadspace=.5em,
    headfont=\normalfont\bfseries,
    headpunct={},
    headformat={\NUMBER.\@\ \NAME.\@\NOTE},
    notefont=\normalfont\bfseries\boldmath,
    notebraces={}{.},
]{definition}
\declaretheoremstyle[
    spaceabove=\parasep, spacebelow=\parasep,
    postheadspace=.5em,
    headfont=\normalfont\bfseries,
    headpunct={},
    headformat={\NUMBER.\@\NOTE},
    notefont=\normalfont\bfseries\boldmath,
    notebraces={}{.},
]{para}

\renewenvironment{proof}[1][\proofname]{\par
    \pushQED{\qed}%
    \normalfont\trivlist
    \item[\hskip\labelsep\bfseries #1\@addpunct{.}]\ignorespaces
}{%
    \popQED\endtrivlist\@endpefalse
}
\makeatother

\declaretheorem[sibling=paragraph, style=para, refname={\S,\S\S}]{para}
\declaretheorem[sibling=paragraph, style=theorem, name=Theorem]{theorem}
\declaretheorem[sibling=paragraph, style=theorem, name=Lemma]{lemma}

\declaretheorem[numbered=no, style=theorem*, name=Theorem]{theorem*}
\declaretheorem[numbered=no, style=theorem*, name=Lemma]{lemma*}

\numberwithin{equation}{paragraph}

\crefdefaultlabelformat{#2{\upshape#1}#3}

\crefformat{equation}{#2{\upshape(#1)}#3}
\crefrangeformat{equation}{{\upshape#3(#1)#4--#5(#2)#6}}
\crefmultiformat{equation}{#2{\upshape(#1)}#3}{ and #2{\upshape(#1)}#3}{, #2{\upshape(#1)}#3}{, and~#2{\upshape(#1)}#3}

\crefformat{enumi}{#2{\upshape#1}#3}
\crefrangeformat{enumi}{{\upshape#3#1#4--#5#2#6}}
\crefmultiformat{enumi}{#2{\upshape#1}#3}{ and #2{\upshape#1}#3}{, #2{\upshape#1}#3}{, and~#2{\upshape#1}#3}

\crefformat{section}{#2{\upshape\S#1}#3}
\crefrangeformat{section}{{\upshape#3\S\S#1#4--#5#2#6}}
\crefmultiformat{section}{{\upshape#2\S\S#1#3}}{ and {\upshape#2#1#3}}{, {\upshape#2#1#3}}{, and~{\upshape#2#1#3}}
\crefformat{subsection}{#2{\upshape\S#1}#3}
\crefrangeformat{subsection}{{\upshape#3\S\S#1#4--#5#2#6}}
\crefmultiformat{subsection}{{\upshape#2\S\S#1#3}}{ and {\upshape#2#1#3}}{, {\upshape#2#1#3}}{, and~{\upshape#2#1#3}}
\crefformat{paragraph}{#2{\upshape\S#1}#3}
\crefrangeformat{paragraph}{{\upshape#3\S\S#1#4--#5#2#6}}
\crefmultiformat{paragraph}{{\upshape#2\S\S#1#3}}{ and {\upshape#2#1#3}}{, {\upshape#2#1#3}}{, and~{\upshape#2#1#3}}
\crefformat{para}{#2{\upshape\S#1}#3}
\crefrangeformat{para}{{\upshape#3\S\S#1#4--#5#2#6}}
\crefmultiformat{para}{{\upshape#2\S\S#1#3}}{ and {\upshape#2#1#3}}{, {\upshape#2#1#3}}{, and~{\upshape#2#1#3}}

\usepackage{enumitem}
\setlist{noitemsep}
\setlist[enumerate]{label=\textnormal{(\roman*)}}

\usepackage[labelsep=period, labelfont={bf}]{caption}
\numberwithin{figure}{section}
\numberwithin{table}{section}

\usepackage[
    backend=biber,
    style=oxnum,
    giveninits,
    maxbibnames=5,
    scnames,
    sorting=nyt
]{biblatex}

\DeclareFieldFormat[article]{version}{\IfDecimal{#1}{version~}{}#1}

%% file: _commands.tex
\newcommand{\bbA}{\mathbb{A}}
\newcommand{\bbC}{\mathbb{C}}
\newcommand{\bbK}{\mathbb{K}}
\newcommand{\bbL}{\mathbb{L}}

\newcommand{\bbP}{\mathbb{P}}
\newcommand{\bbQ}{\mathbb{Q}}
\newcommand{\bbR}{\mathbb{R}}
\newcommand{\bbT}{\mathbb{T}}
\newcommand{\bbZ}{\mathbb{Z}}
\newcommand{\calA}{\mathcal{A}}

\newcommand{\calD}{\mathcal{D}}
\newcommand{\calE}{\mathcal{E}}

\newcommand{\calI}{\mathcal{I}}

\newcommand{\calO}{\mathcal{O}}
\newcommand{\calP}{\mathcal{P}}

\newcommand{\calU}{\mathcal{U}}

\newcommand{\calX}{\mathcal{X}}
\newcommand{\calY}{\mathcal{Y}}
\newcommand{\calZ}{\mathcal{Z}}
\newcommand{\cat}[1]{\smash[b]{\textsf{\textup{#1}}}}
\newcommand{\CF}{\mathrm{CF}}

\newcommand{\colim}{\operatornamewithlimits{colim}}
\newcommand{\Crit}{\mathrm{Crit}}
\newcommand{\cwP}{{\smash{^\mathrm{cw}}\mathbb{P}}}

\newcommand{\dFilt}{\mathrm{dFilt}}
\newcommand{\dGrad}{\mathrm{dGrad}}
\newcommand{\dMap}{\mathrm{dMap}}
\newcommand{\DT}{\mathrm{DT}}

\newcommand{\ev}{\mathrm{ev}}

\newcommand{\Filt}{\mathrm{Filt}}
\newcommand{\frg}{\mathfrak{g}}
\newcommand{\frM}{\mathfrak{M}}
\newcommand{\frp}{\mathfrak{p}}

\newcommand{\frL}{\mathfrak{L}}
\newcommand{\frX}{\mathfrak{X}}
\newcommand{\frY}{\mathfrak{Y}}
\newcommand{\GL}{\mathrm{GL}}
\newcommand{\Ga}{\mathbb{G}_{\mathrm{a}}}
\newcommand{\Gm}{\mathbb{G}_{\mathrm{m}}}
\newcommand{\gr}{\mathrm{gr}}
\newcommand{\Grad}{\mathrm{Grad}}

\newcommand{\id}{\mathrm{id}}

\newcommand{\Kvar}{K_{\mathrm{var}}}
\newcommand{\Kvarm}{K_{\mathrm{var}}^{\smash{\hat{\upmu}}}}
\newcommand{\longhookrightarrow}{\lhook\joinrel\longrightarrow}

\newcommand{\longrightrightarrows}{\mathrel{\raisebox{-.4ex}{$\overset{\smash[b]{\raisebox{-.7ex}[.2ex]{$\longrightarrow$}}}{\longrightarrow}$}}}
\newcommand{\longsimto}{\mathrel{\overset{\smash{\raisebox{-.8ex}{$\sim$}}\mspace{3mu}}{\longrightarrow}}}
\newcommand{\Map}{\mathrm{Map}}

\newcommand{\MF}{\mathrm{MF}}
\newcommand{\Mhat}{\hat{\mathrm{M}}}
\newcommand{\Mm}{\mathrm{M}^{\smash{\hat{\upmu}}}}
\newcommand{\Mmhat}{\hat{\mathrm{M}}{}^{\smash{\hat{\upmu}}}}
\newcommand{\mot}{{\smash{\mathrm{mot}}}}
\newcommand{\muhat}{{\smash{\hat{\upmu}}}}
\newcommand{\numberthis}{\addtocounter{equation}{1}\tag{\theequation}}
\newcommand{\oline}[1]{{\smash{\overline{\smash{#1}\vphantom{X}}}}}
\newcommand{\op}{{\smash{\mathrm{op}}}}\newcommand{\Perf}{\cat{Perf}}

\newcommand{\pn}[1]{{\smash{(#1)}}}
\newcommand{\pr}{\operatorname{\smash[b]{\mathrm{pr}}}}
\renewcommand{\qed}{\hfill\ensuremath{\square}}
\newcommand{\rank}{\operatorname{rank}}
\newcommand{\red}{{\smash{\mathrm{red}}}}
\newcommand{\ssf}{\mathrm{sf}}

\newcommand{\simto}{\mathrel{\overset{\smash{\raisebox{-.8ex}{$\sim$}}\mspace{3mu}}{\to}}}

\newcommand{\slsl}{\mathord{\mathrlap{/}\mspace{3mu}/}}

\newcommand{\Spec}{\operatorname{Spec}}

\newcommand{\tot}{\mathrm{tot}}

\newcommand{\upH}{\mathrm{H}}

\newcommand{\upM}{\mathrm{M}}

\newcommand{\vdim}{\operatorname{vdim}}
\newcommand{\wBl}{{\smash{^\mathrm{w}}\mathrm{Bl}}}
\newcommand{\wP}{{\smash{^\mathrm{w}}\mathbb{P}}}

\renewcommand{\leq}{\leqslant}
\renewcommand{\geq}{\geqslant}

\newcommand{\titleminusone}{(\hspace{.02em}\scalebox{1.1}{\textminus}1\hspace{.02em})}

%% file: abstract.tex
We prove a motivic integral identity
relating the motivic Behrend function
of a $(-1)$\=/shifted symplectic stack
to that of its stack of graded points.
This generalizes analogous identities
for moduli stacks of objects in
$3$\=/Calabi--Yau abelian categories
obtained by Kontsevich--Soibelman
and Joyce--Song,
which are crucial in proving wall-crossing formulae
for Donaldson--Thomas invariants.
We expect our identity to be useful
in extending motivic Donaldson--Thomas theory
to general $(-1)$\=/shifted symplectic stacks.

%% file: intro.tex
\addtocounter{subsection}{1}

\begin{para}
    Let~$\bbK$ be an algebraically closed field of characteristic zero,
    and let~$\frX$ be a \emph{$(-1)$\=/shifted symplectic derived algebraic stack} over~$\bbK$,
    in the sense of Pantev--To\"en--Vaqui\'e--Vezzosi~\cite{PTVV2013},
    such that its classical truncation is a classical algebraic $1$-stack.

    One of the motivating examples of such a stack
    is the moduli stack $\frX = \frM_{\calA}$
    of objects in a $\bbK$\=/linear $3$\=/Calabi--Yau category~$\calA$,
    such as the category of coherent sheaves
    on a smooth projective Calabi--Yau threefold.
    Such stacks are of great interest in Donaldson--Thomas theory,
    studied by Joyce--Song~\cite{JoyceSong2012},
    Kontsevich--Soibelman~\cite{KontsevichSoibelman2008},
    and many others.

    Starting from~$\frX$,
    one can consider the derived mapping stacks%
    \begin{align*}
        \Grad (\frX) & = \Map ([* / \Gm], \frX) \ , \\
        \Filt (\frX) & = \Map ([\bbA^1 / \Gm], \frX) \ ,
    \end{align*}
    called the \emph{stack of graded points}
    and the \emph{stack of filtered points} of~$\frX$, respectively,
    as in Halpern-Leistner~\cite{HalpernLeistnerInstability}.
    For example, if $\frX = \frM_{\calA}$ as above,
    then $\Grad (\frX)$ is the moduli stack of $\bbZ$\=/graded objects in~$\calA$,
    and $\Filt (\frX)$ is the moduli stack of $\bbZ$\=/filtered objects in~$\calA$.
    In particular, there are inclusions%
    \begin{align}
        \label{eq-intro-grad-eg}
        \frM_{\calA} \times \frM_{\calA} & \longhookrightarrow \Grad (\frM_{\calA}) \ , \\
        \label{eq-intro-filt-eg}
        \frM_{\mathrm{Exact} (\calA)} & \longhookrightarrow \Filt (\frM_{\calA})
    \end{align}
    as open and closed substacks,
    i.e.~disjoint unions of connected components,
    where $\frM_{\mathrm{Exact} (\calA)}$ is the moduli stack of
    short exact sequences in~$\calA$.
    These substacks can be given by,
    for example, graded and filtered objects
    that only have non-trivial factors
    in degrees~$0$ and~$1$.
\end{para}

\begin{para}
    We show in \cref{thm-filt-lag} that
    there is a $(-1)$\=/shifted Lagrangian correspondence%
    \begin{equation}
        \label{eq-intro-filt-corr}
        \smash{ \Grad (\frX) \overset{\gr}{\longleftarrow}
        \Filt (\frX) \overset{\ev_1}{\longrightarrow} \frX \ . }
    \end{equation}
    When $\frX = \frM_{\calA}$,
    the morphisms~$\gr$ and~$\ev_1$ send a filtered object to
    its associated graded object and its total object, respectively.
    In particular, restricting it to the substacks
    \crefrange{eq-intro-grad-eg}{eq-intro-filt-eg},
    this gives the $(-1)$\=/shifted Lagrangian correspondence%
    \begin{equation}
        \label{eq-intro-filt-corr-eg}
        \smash{ \frM_{\calA} \times \frM_{\calA} \overset{(p_1, \, p_3)}{\longleftarrow}
        \frM_{\mathrm{Exact} (\calA)} \overset{p_2}{\longrightarrow} \frM_{\calA} }
    \end{equation}
    as in Brav--Dyckerhoff~\cite{BravDyckerhoff2021},
    where $p_1, p_2, p_3 \colon \frM_{\mathrm{Exact} (\calA)} \to \frM_{\calA}$
    send a short exact sequence to its three respective terms.

    The correspondence~\cref{eq-intro-filt-corr-eg}
    has proved to be useful in enumerative geometry.
    It lies in the heart of the construction of
    Hall-algebra-type algebraic structures,
    including \emph{motivic Hall algebras}
    studied by Joyce~\cite{Joyce2007II},
    \emph{cohomological Hall algebras}
    introduced by Kontsevich--Soibelman~\cite{KontsevichSoibelman2011},
    and \emph{Joyce vertex algebras}
    constructed by Joyce~\cite{Joyce2021,JoyceHall}.
    These structures are closely related to
    Donaldson--Thomas invariants and other enumerative invariants.
\end{para}

\begin{para}
    \label{para-intro-behrend}
    Following a series of works
    \cite{BBJ2019,Joyce2015,BJM2019,BBBBJ2015,BBDJS2015}
    by Joyce and his collaborators,
    it is known that a $(-1)$\=/shifted symplectic stack~$\frX$
    can be locally modelled as
    \emph{derived critical loci} of functions on smooth stacks.
    When $\frX$ is equipped with \emph{orientation data},
    one can define an element~$\nu^\mot_{\frX}$ in the
    ring of monodromic motives on~$\frX$,
    which we call the \emph{motivic Behrend function} of~$\frX$.
    It is locally modelled by the \emph{motivic vanishing cycle}
    defined by Denef--Loeser~\cite{DenefLoeser2001},
    and is a motivic enhancement of the \emph{Behrend function}
    $\nu_{\frX} \colon \frX \to \bbZ$
    introduced by Behrend~\cite{Behrend2009}
    and extended by Joyce--Song~\cite[\S4.1]{JoyceSong2012} to algebraic stacks.

    When $\frX = \frM_{\calA}$ as above, the element $\nu^\mot_{\frM_{\calA}}$
    was considered by
    Kontsevich--Soibelman~\cite{KontsevichSoibelman2008},
    and is important in the Donaldson--Thomas theory of~$\calA$.
    Given a stability condition on~$\calA$,
    if $\frM_\alpha \subset \frM_{\calA}$ is a component,
    then the \emph{motivic Donaldson--Thomas invariant}
    of the class~$\alpha$ is the monodromic motive
    given by the \emph{motivic integral}%
    \begin{equation}
        \DT^\mot_\alpha =
        \int_{\frM_\alpha} {} (\bbL^{1/2} - \bbL^{-1/2}) \cdot \epsilon_\alpha \cdot \nu^\mot_{\frM_{\calA}} \ ,
    \end{equation}
    where $\epsilon_\alpha$ is a \emph{weight function}
    encoding the data of the stability condition.
    Similarly, the numerical \emph{Donaldson--Thomas invariant}~$\DT_\alpha$
    is given by a weighted Euler characteristic%
    \begin{equation}
        \DT_\alpha = \int_{\frM_\alpha} {} (\bbL^{1/2} - \bbL^{-1/2}) \cdot \epsilon_\alpha \cdot \nu_{\frM_{\calA}} \, d \chi \ .
    \end{equation}
    See Joyce--Song~\cite{JoyceSong2012} for more details.
\end{para}

\begin{para}
    The main result of this paper, \cref{thm-behrend-main},
    states that the motivic Behrend functions of~$\frX$ and~$\Grad (\frX)$
    are related via the correspondence~\cref{eq-intro-filt-corr},
    by the identity%
    \begin{equation}
        \label{eq-intro-behrend-main}
        \gr_! \circ \ev_1^* (\nu^\mot_{\frX})
        = \bbL^{\vdim \Filt (\frX) / 2} \cdot
        \nu^\mot_{\Grad (\frX)} \ ,
    \end{equation}
    as an identity of monodromic motives on~$\Grad (\frX)$,
    where $\vdim \Filt (\frX)$ is the virtual dimension of~$\Filt (\frX)$.
    In particular, evaluating this at a graded point~$\gamma \in \Grad (\frX)$
    gives the motivic integral identity%
    \begin{equation}
        \label{eq-intro-behrend-pointwise}
        \int \limits_{\varphi \in \gr^{-1} (\gamma)} \nu^\mot_{\frX} (\ev_1 (\varphi))
        = \bbL^{\vdim \Filt (\frX) / 2} \cdot \nu^\mot_{\Grad (\frX)} (\gamma) \ ,
    \end{equation}
    as an identity of monodromic motives over~$\bbK$.

    We prove this identity by first
    proving a local version of it in \cref{thm-behrend-local},
    which is, roughly speaking, the special case
    of~\cref{eq-intro-behrend-main} when $\frX = [\Crit (f) / \Gm]$
    is a derived critical locus,
    where $f$~is a $\Gm$-invariant function
    on a $\Gm$-equivariant smooth $\bbK$\=/variety.
    Our proof of this local version involves the theory of
    \emph{nearby and vanishing cycles}
    for rings of motives on algebraic stacks,
    which we develop in \cref{subsec-nearby-stacks},
    with \cref{thm-nearby-properties} as a main result.
    Then, we prove the global version of the identity
    by gluing together the local models.
\end{para}

\begin{para}
    We explain relations between the identity~\cref{eq-intro-behrend-main}
    and known results and conjectures in the literature.

    Firstly, this identity can be seen as
    a global version and a generalization
    of an integral identity conjectured by
    Kontsevich--Soibelman~\cite[Conjecture~4]{KontsevichSoibelman2008},
    later proved by L\^e~\cite{Le2015}.
    The local version of our identity, \cref{thm-behrend-local},
    is stated in a form similar to Kontsevich--Soibelman's identity,
    and generalizes it by removing the assumption that
    the torus action only has weights~$-1$, $0$, and~$1$.
    It is crucial that this assumption is removed
    in order for the identity to serve as
    a local model for~\cref{eq-intro-behrend-main}
    for general $(-1)$\=/shifted symplectic stacks,
    not only for stacks of the form~$\frM_{\calA}$.

    Kontsevich--Soibelman then used their identity
    to prove \cite[Theorem~8]{KontsevichSoibelman2008},
    which can be seen as a special case
    of our identity~\cref{eq-intro-behrend-main} when $\frX = \frM_{\calA}$ as above.
    Their theorem is a key ingredient in proving
    \emph{wall-crossing formulae} of motivic Donaldson--Thomas invariants,
    governing the behaviour of these invariants
    under changes of stability conditions.

    Secondly, by taking the Euler characteristic of our identity,
    we obtain numerical integral identities in \cref{thm-behrend-num}.
    These identities are direct generalizations of
    the Behrend function identities of
    Joyce--Song~\cite[Theorem~5.11]{JoyceSong2012}
    to general $(-1)$\=/shifted symplectic stacks.

    Thirdly, the identity~\cref{eq-intro-behrend-main}
    is related to a conjecture on perverse sheaves,
    sometimes known as the \emph{Joyce conjecture},
    formulated in Joyce--Safronov~%
    \cite[Conjecture~1.1]{JoyceSafronov2019}
    and in Amorim--Ben-Bassat~%
    \cite[\S5.3]{AmorimBenBassat2017}.
    One form of the conjecture states that
    for an oriented $(-1)$\=/shifted Lagrangian correspondence%
    \begin{equation}
        \label{eq-intro-joyce-corr}
        \frX \overset{f}{\longleftarrow}
        \frL \overset{g}{\longrightarrow} \frY \ ,
    \end{equation}
    under certain assumptions,
    there should exist a natural morphism%
    \begin{equation}
        \label{eq-intro-joyce-conj}
        \mu_{\frL} \colon f_! \circ g^* (\calP_{\frY})
        \longrightarrow \calP_{\frX} [-{\vdim \frL}] \ ,
    \end{equation}
    satisfying certain properties,
    where $\calP_{\frX}$ and $\calP_{\frY}$ are the perverse sheaves
    constructed in Ben-Bassat--Brav--Bussi--Joyce~%
    \cite[Theorem~4.8]{BBBBJ2015},
    sometimes called the \emph{Donaldson--Thomas perverse sheaves}.
    They can be seen as analogues of the motivic Behrend functions
    $\nu^\mot_{\frX}$ and $\nu^\mot_{\frY}$
    in cohomological Donaldson--Thomas theory.

    In the special case when the correspondence~\cref{eq-intro-joyce-corr}
    is taken to be the correspondence~\cref{eq-intro-filt-corr},
    a recent result of Kinjo--Park--Safronov~\cite[Theorem~B]{KinjoParkSafronov}
    shows that there is a natural isomorphism
    of the form~\cref{eq-intro-joyce-conj},
    strengthening the Joyce conjecture in this special case.
    In this sense, the identity~\cref{eq-intro-behrend-main}
    can be seen as a motivic analogue
    of this version of the Joyce conjecture.
\end{para}

\begin{para}
    In a series of works in progress
    \cite{BuHalpernLeistnerIbanezNunezKinjo,BuIbanezNunezKinjoII,BuIbanezNunezKinjoIII},
    the author and his collaborators
    plan to extend the definition
    of the weight functions~$\epsilon_\alpha$
    mentioned in~\cref{para-intro-behrend}
    to general algebraic stacks,
    which will enable us to define motivic Donaldson--Thomas invariants
    for general $(-1)$\=/shifted symplectic stacks
    equipped with extra data similar to a stability condition,
    under mild assumptions.

    The integral identity~\cref{eq-intro-behrend-main}
    will be crucial in proving \emph{wall-crossing formulae}
    for these general Donaldson--Thomas invariants,
    which relate the invariants
    for the same stack under different stability conditions.
    Wall-crossing formulae have seen many important applications
    in enumerative geometry,
    since they impose a very strong constraint on the structure
    of enumerative invariants,
    and can sometimes be used to compute them explicitly.
    See Joyce--Song~\cite{JoyceSong2012} and Kontsevich--Soibelman~\cite{KontsevichSoibelman2008}
    in the case of motivic Donaldson--Thomas theory,
    and \mbox{\cite{BojkoLimMoreira2024,Bu2023Curves,GrossJoyceTanaka2022,Joyce2021},}
    etc., for applications in other contexts.
    Using wall-crossing formulae for general stacks,
    we hope to generalize many of the applications mentioned above
    to general $(-1)$\=/shifted symplectic stacks.

    For example, as a special case of this generalized theory,
    the author~\cite{Bu2023I}
    defines the weight functions for the moduli stack of
    \emph{self-dual objects} in certain \emph{self-dual $\bbK$\=/linear categories}.
    Such stacks include the moduli stack of
    principal orthogonal or symplectic bundles
    on a smooth projective variety,
    or a certain compactification of it.
    Combined with the contents of the present work,
    it will become possible to define and study
    motivic Donaldson--Thomas invariants
    for type~B/C/D structure groups on a Calabi--Yau threefold,
    and prove wall-crossing formulae for them.
    The author plans to report on this in a future paper.
\end{para}

\begin{para}[Acknowledgements]
    The author would like to thank
    Andrés Ibáñez Núñez,
    Dominic Joyce,
    and Tasuki Kinjo,
    for helpful discussions and comments.
    The author also thanks the anonymous referee
    for their valuable suggestions.

    The author is grateful to the
    Mathematical Institute, University of Oxford,
    for its support during the preparation of this paper.
\end{para}

\begin{para}[Conventions]
    \label{para-conventions}
    Throughout this paper,
    we will use the following notations, terminology, and conventions.

    \begin{itemize}
        \item
            $\bbK$ is an algebraically closed field
            of characteristic zero.
        \item
            A \emph{$\bbK$\=/variety}
            is a separated $\bbK$\=/scheme of finite type.
        \item
            A \emph{reductive group} over~$\bbK$
            is a linear algebraic group over~$\bbK$
            that is linearly reductive, and is allowed to be disconnected.
        \item
            All \emph{$\bbK$\=/schemes}, \emph{algebraic spaces} over~$\bbK$,
            and \emph{algebraic stacks} over~$\bbK$
            are assumed to be quasi-separated and locally of finite type.
            Algebraic stacks are assumed
            to have separated diagonal.
        \item
            A \emph{derived algebraic stack} over~$\bbK$
            is a derived stack over~$\bbK$ that has an open cover by
            \emph{geometric stacks}
            in the sense of To\"en--Vezzosi~%
            \cite[\S1.3.3]{ToenVezzosi2008},
            and is assumed to be locally almost of finite presentation.
            When it is locally of finite presentation,
            its cotangent complex is perfect,
            and the rank of its cotangent complex is called its
            \emph{virtual dimension}.
        \item
            An \emph{$s$\=/shifted symplectic stack} over~$\bbK$,
            where $s \in \bbZ$,
            is a derived algebraic stack
            locally of finite presentation over~$\bbK$,
            equipped with an \emph{$s$\=/shifted symplectic structure}
            in the sense of Pantev--To\"en--Vaqui\'e--Vezzosi~%
            \cite[\S1]{PTVV2013}.
    \end{itemize}
\end{para}

%% file: mot.tex
The main purpose of this section is to define and
study the \emph{motivic nearby and vanishing cycle maps}
on certain rings of motives over algebraic stacks.
These are based on the construction of the \emph{motivic Milnor fibre}
by Denef--Loeser~\cite{DenefLoeser1998,DenefLoeser2001,DenefLoeser2002}
and Looijenga~\cite{Looijenga2002},
and are a generalization of the work of Bittner~\cite{Bittner2005}
from the case of varieties to that of stacks.
We define these maps in \cref{thm-nearby-stacks},
and prove a useful property in \cref{thm-nearby-properties},
which will be used in the proof of our main results.

In \cref{subsec-nu-mot},
we discuss the \emph{motivic Behrend function} introduced by
Bussi--Joyce--Meinhardt~\cite{BJM2019}
and Ben-Bassat--Brav--Bussi--Joyce~\cite{BBBBJ2015},
and slightly generalize their construction
by weakening the assumptions on the stack.

\subsection{Rings of motives}
\label{subsec-mot-bg}

\begin{para}
    We provide background on rings of motives
    over schemes, algebraic spaces, and algebraic stacks.
    The case of stacks was first due to Joyce~\cite{Joyce2007Stack},
    where these rings of motives, along with several variations,
    were constructed and called rings of \emph{stack functions}.

    To avoid repetition,
    we only state the majority of definitions and results
    for stacks, with the understanding that schemes and algebraic spaces
    are special cases of algebraic stacks.
\end{para}

\begin{para}[Stacks with affine stabilizers]
    Let $\calX$ be an algebraic stack over~$\bbK$.
    We say that $\calX$ has \emph{affine stabilizers},
    if for any field~$k$ and any point $x \in \calX (k)$,
    the stabilizer group of~$\calX$ at~$x$
    is an affine algebraic group over~$k$.
\end{para}

\begin{para}[Rings of motives]
    \label{para-def-kvar}
    Let $\calX$ be an algebraic stack over~$\bbK$
    with affine stabilizers.
    Define the \emph{Grothendieck ring of varieties} over $\calX$
    to be the abelian group
    \begin{equation}
        \Kvar (\calX) = \mathop{\hat{\bigoplus}}_{Z \to \calX} {} \bbZ \cdot [Z] \, \Big/ {\sim} \ ,
    \end{equation}
    where we run through all morphisms $Z \to \calX$
    with $Z$ a $\bbK$\=/variety,
    and $\hat{\oplus}$ means we take the set of \emph{locally finite sums},
    that is, possibly infinite sums $\sum_{Z \to \calX} n_Z \cdot [Z]$,
    such that for each quasi-compact open substack $\calU \subset \calX$,
    there are only finitely many $Z$ such that
    $n_Z \neq 0$ and $Z \times_{\calX} \calU \neq \varnothing$.
    The relation $\sim$ is generated by
    $[Z] \sim [Z'] + [Z \setminus Z']$ for closed subschemes $Z' \subset Z$.

    One can define multiplication on $\Kvar (\calX)$
    by taking the fibre product over~$\calX$,
    making it into a commutative ring,
    possibly non-unital when $\calX$ is not an algebraic space.
    It is also a (possibly non-unital) commutative $\Kvar (\bbK)$-algebra,
    where $\Kvar (\bbK) = \Kvar (\Spec \bbK)$,
    with the action given by the product.

    Let $\bbL = [\bbA^1] \in \Kvar (\bbK)$,
    and define \emph{rings of motives} over $\calX$,
    \begin{align}
        \upM (\calX) & = \Kvar (\calX) \underset{\Kvar (\bbK)}{\mathbin{\hat{\otimes}}}
        \Kvar (\bbK) \, [\bbL^{-1}] \big/
        (\bbL - 1) \text{-torsion} \ ,
        \\
        \label{eq-def-mhat}
        \Mhat (\calX) & = \Kvar (\calX) \underset{\Kvar (\bbK)}{\mathbin{\hat{\otimes}}}
        \Kvar (\bbK) \, [\bbL^{-1}, (\bbL^k - 1)^{-1}] \ ,
    \end{align}
    where $\hat{\otimes}$ means we take the set of locally finite sums,
    that is, possibly infinite sums $\sum_{Z \to \calX} [Z] \otimes a_Z$,
    such that the family of~$Z$ with~$a_Z \neq 0$ is locally finite on~$\calX$,
    and in~\cref{eq-def-mhat}, we invert $\bbL^k - 1$ for all $k \geq 1$.
\end{para}

\begin{para}[Motives of algebraic spaces and stacks]
    \label{para-mot-of-stacks}
    Let $\calX$ be as above.
    For an algebraic space~$Z$ and a morphism~$Z \to \calX$ of finite type,
    one can also assign a class $[Z] \in \Kvar (\calX)$,
    extending the usual definition for varieties,
    since $Z$ can be stratified by varieties.

    Furthermore, as in Joyce~\cite{Joyce2007Stack}
    or Ben-Bassat--Brav--Bussi--Joyce~\cite[\S5.3]{BBBBJ2015},
    for any finite type morphism of algebraic stacks $\calZ \to \calX$,
    where $\calZ$ has affine stabilizers,
    one can assign a class $[\calZ] \in \Mhat (\calX)$,
    which agrees with the usual one when $\calZ$ is a variety,
    and satisfies the relation $[\calZ] = [\calZ'] + [\calZ \setminus \calZ']$
    for closed substacks $\calZ' \subset \calZ$.
    In particular, the class $[\calX] \in \Mhat (\calX)$
    is the multiplicative unit of the ring~$\Mhat (\calX)$.

    For an algebraic stack $\calX$ over~$\bbK$ of finite type,
    with affine stabilizers, one thus have a class
    $[\calX] \in \Mhat (\bbK)$, called the \emph{motive} of~$\calX$.
\end{para}

\begin{para}[Pullbacks and pushforwards]
    Let $\calX, \calY$ be algebraic stacks over~$\bbK$ with affine stabilizers,
    and let $f \colon \calX \to \calY$ be a morphism.

    There is a pullback map
    \begin{equation*}
        f^* \colon \Mhat (\calY) \longrightarrow \Mhat (\calX) \ ,
    \end{equation*}
    which is an $\Mhat (\bbK)$-algebra homomorphism,
    given on generators by $f^* [Z] = [Z \times_{\calY} \calX]$,
    where the right-hand side is defined as in~\cref{para-mot-of-stacks}.
    Pulling back respects composition of morphisms.
    If, moreover, $f$ is representable,
    then there are pullback maps
    \begin{align*}
        f^* \colon \Kvar (\calY) & \longrightarrow \Kvar (\calX) \ ,
        \\
        f^* \colon \upM (\calY) & \longrightarrow \upM (\calX) \ ,
    \end{align*}
    which are $\Kvar (\bbK)$- and $\upM (\bbK)$-algebra homomorphisms,
    respectively, and are defined similarly.

    On the other hand, if $f$ is of finite type,
    then there are pushforward maps
    \begin{align*}
        f_! \colon \Kvar (\calX) & \longrightarrow \Kvar (\calY) \ ,
        \\
        f_! \colon \upM (\calX) & \longrightarrow \upM (\calY) \ ,
        \\
        f_! \colon \smash{\Mhat (\calX)} & \longrightarrow \smash{\Mhat (\calY)} \ ,
    \end{align*}
    which are $\Kvar (\bbK)$-, $\upM (\bbK)$-,
    and $\Mhat (\bbK)$-module homomorphisms, respectively,
    given on generators by $f_! [Z] = [Z]$.
    Pushing forward respects composition of morphisms.

    In particular, when $\calX$ is of finite type,
    pushing forward along the structure morphism
    $\calX \to \Spec \bbK$ is sometimes called \emph{motivic integration},
    and denoted by
    \begin{equation*}
        \int_{\calX} {} (-) \colon
        \Mhat (\calX) \longrightarrow \Mhat (\bbK) \ .
    \end{equation*}
\end{para}

\begin{para}[Base change and projection formulae]
    \label{para-mot-bc-pf}
    Suppose we have a pullback diagram
    \begin{equation*}
        \begin{tikzcd}
            \calX' \ar[r, "f'"] \ar[d, "g'"']
            \ar[dr, phantom, pos=.2, "\ulcorner"] &
            \calY' \ar[d, "g"]
            \\
            \calX \ar[r, "f"] &
            \calY \rlap{ ,}
        \end{tikzcd}
    \end{equation*}
    where $\calX, \calY, \calX', \calY'$ are algebraic stacks over~$\bbK$
    with affine stabilizers,
    and $f$ is of finite type.
    Then we have the \emph{base change formula}
    \begin{equation}
        \label{eq-mot-bc}
        g^* \circ f_! = f'_! \circ g'^*
    \end{equation}
    on $\Mhat (-)$, as in \cite[Theorem~3.5]{Joyce2007Stack}.
    Moreover, if $g$ is representable,
    then this also holds for $\Kvar (-)$ and $\upM (-)$.

    Let $f \colon \calX \to \calY$ be as above.
    We have the \emph{projection formula}
    \begin{equation}
        \label{eq-mot-pf}
        f_! (a \cdot f^* (b)) = f_! (a) \cdot b
    \end{equation}
    for all $a \in \Mhat (\calX)$ and $b \in \Mhat (\calY)$,
    which can be verified directly on generators.
    Moreover, if~$f$ is representable,
    then this also holds for $\Kvar (-)$ and $\upM (-)$.
\end{para}

\begin{para}[Motives of principal bundles]
    \label{para-pr-bun-rel}
    Following Serre~\cite[\S4]{Serre1958},
    an algebraic group $G$ over~$\bbK$ is \emph{special}
    if all principal $G$-bundles over a $\bbK$\=/scheme
    are Zariski locally trivial.
    For example, the groups $\GL (n)$ and $\Ga$ are special;
    semidirect products of special groups are special;
    disconnected groups are not special.
    See also Joyce~\cite[Definition~2.1]{Joyce2007Stack} for related discussions.

    For a special group~$G$ and a principal $G$-bundle
    $\pi \colon \calY \to \calX$, where $\calX$ is an algebraic stack over~$\bbK$
    with affine stabilizers, we have the relation
    \begin{equation}
        \label{eq-pr-bun-rel}
        \pi_! \circ \pi^* = [G] \cdot \id
    \end{equation}
    in~$\Kvar (\calX)$, $\upM (\calX)$, and $\Mhat (\calX)$,
    which can be deduced from the definition of special groups.
    In particular, we have $[\calY] = [G] \cdot [\calX]$ in~$\Mhat (\calX)$.

    Note, however, that such a principal $G$-bundle~$\pi$ itself
    is not necessarily Zariski locally trivial,
    but it becomes so after a base change to any scheme.

    The relation~\cref{eq-pr-bun-rel}
    is not necessary true if $G$ is not special.
    For example, consider the principal $\bbZ_2$-bundle
    $\Gm \to \Gm$ given by $t \mapsto t^2$.
    Then the equality cannot hold, since
    $[\Gm] \neq 2 \cdot [\Gm]$.
\end{para}

\begin{para}[The Euler characteristic]
    \label{para-euler-char}
    As in Joyce~\cite[Example~6.3]{Joyce2007Stack}, there is a ring map
    \begin{equation*}
        \chi \colon \upM (\bbK) \longrightarrow \bbZ \ ,
    \end{equation*}
    sending each generator $[Z]$ to its Euler characteristic,
    and sending $\bbL$ to $1$.
    This extends naturally to a map
    $\chi \colon \Mhat (\bbK) \to \bbQ \cup \{ \infty \}$,
    sending $1 / (1 + \bbL + \cdots + \bbL^{k - 1})$ to $1/k$
    for each $k \geq 1$, and sending elements not in
    $\upM (\bbK) [(1 + \bbL + \cdots + \bbL^{k - 1})^{-1} : k \geq 1]$ to $\infty$.

    For an algebraic stack~$\calX$ over~$\bbK$ of finite type,
    with affine stabilizers, we have a map
    \begin{equation*}
        \int_{\calX} {} (-) \, d \chi \colon
        \Mhat (\calX) \longrightarrow \bbQ \cup \{ \infty \} \ ,
    \end{equation*}
    defined by pushing forward along $\calX \to \Spec \bbK$,
    and then taking the Euler characteristic.
    We have $\int_{\calX} {} a \, d \chi \in \bbZ$ for all $a \in \upM (\calX)$.
\end{para}

\begin{para}[Constructible functions]
    \label{para-cf}
    For an algebraic stack~$\calX$ over~$\bbK$,
    a \emph{constructible function} on $\calX$
    is a map of sets
    \begin{equation*}
        a \colon | \calX | \longrightarrow \bbZ \ ,
    \end{equation*}
    where $| \calX |$ is the underlying topological space of~$\calX$,
    such that for any $c \in \bbZ$,
    the preimage $a^{-1} (c)$ is a locally constructible subset of~$| \calX |$.
    The abelian group of constructible functions on~$\calX$
    is denoted by $\CF (\calX)$.

    There is an Euler characteristic map
    \begin{equation*}
        \chi \colon \tilde{\upM} (\calX) \longrightarrow \CF (\calX) \ ,
    \end{equation*}
    where $\tilde{\upM} (\calX) \subset \Mhat (\calX)$
    is the $\upM (\bbK)$-subalgebra generated by classes
    $[\calZ]$ for representable morphisms $\calZ \to \calX$,
    and $\chi$ is given by taking the fibrewise Euler characteristic,
    as in Joyce~\cite[Definition~3.2]{Joyce2007Stack},
    where this map was denoted by $\pi_{\calX}^{\smash{\mathrm{stk}}}$.

    One can also define pullback and pushforward maps on $\CF (-)$
    for representable morphisms,
    where pushing forward also requires the morphism to be of finite type.
    They are compatible with the map~$\chi$,
    and satisfy the base change and projection formulae
    as in~\cref{para-mot-bc-pf}, for representable morphisms.
\end{para}

\begin{para}[Rings of monodromic motives]
    \label{para-mon-mot}
    Let $\muhat = \lim \upmu_n$ be the projective limit
    of the groups~$\upmu_n$ of roots of unity.
    For a $\bbK$\=/scheme~$Z$, a \emph{good $\muhat$-action}
    on~$Z$ is one that factors through $\upmu_n$ for some~$n$,
    such that each orbit is contained in
    an affine open subscheme.

    Let $\calX$ be an algebraic stack over~$\bbK$ with affine stabilizers.
    Define the \emph{monodromic Grothendieck ring of varieties} over~$\calX$
    to be the abelian group
    \begin{equation}
        \Kvarm (\calX) = \mathop{\hat{\bigoplus}}_{Z \to \calX} {}
        \bbZ \cdot [Z]^\muhat \, \Big/ {\sim} \ ,
    \end{equation}
    where $\hat{\oplus}$ indicates taking locally finite sums
    as in~\cref{para-def-kvar},
    and we sum over all morphisms $Z \to \calX$
    with $Z$ a $\bbK$\=/variety with a good $\muhat$-action
    that is compatible with the trivial $\muhat$-action on~$\calX$.
    The relation $\sim$ is generated by
    $[Z]^\muhat \sim [Z']^\muhat + [Z \setminus Z']^\muhat$
    for $\muhat$-invariant closed subschemes $Z' \subset Z$,
    and $[Z \times V]^\muhat \sim [Z \times \bbA^n]^\muhat$
    for a $\muhat$-representation~$V$ of dimension~$n$,
    where the projections to~$\calX$ factor through~$Z$,
    and $\muhat$~acts on $\bbA^n$ trivially.
    See Looijenga~\cite[\S5]{Looijenga2002}
    and Ben-Bassat--Brav--Bussi--Joyce~\cite[\S5]{BBBBJ2015}.

    Using the $\Kvar (\bbK)$-module structure on~$\Kvarm (\calX)$,
    we define \emph{rings of monodromic motives} over~$\calX$,
    \begin{align}
        \Mm (\calX) & = \Kvarm (\calX) \underset{\Kvar (\bbK)}{\mathbin{\hat{\otimes}}}
        \Kvar (\bbK) \, [\bbL^{-1}] \big/
        (\bbL - 1) \text{-torsion} \ ,
        \\
        \Mmhat (\calX) & = \Kvarm (\calX) \underset{\Kvar (\bbK)}{\mathbin{\hat{\otimes}}}
        \Kvar (\bbK) \, [\bbL^{-1}, (\bbL^k - 1)^{-1}]
        \Big/ {\approx} \ ,
    \end{align}
    where $\hat{\otimes}$ indicates taking locally finite sums
    as in~\cref{para-def-kvar}.
    The ring structures and the relation $\approx$ are defined below.

    We consider multiplication on $\Kvarm (\calX)$
    denoted by `$\odot$' in \cite[Definition~5.3]{BBBBJ2015};
    see there for the definition.
    We will denote this by `${} \cdot {}$'.
    This makes $\Kvarm (\calX)$ into a ring,
    possibly non-unital when $\calX$ is not an algebraic space.
    Note that this is \emph{not} given by the fibre product,
    although the latter does define a different ring structure.
    The relation~$\approx$ is defined as in
    \cite[Definitions~5.5 and~5.13]{BBBBJ2015},
    denoted by `$I^{\mathrm{st}, \muhat}_{\smash{\calX}}$' there,
    and is imposed so that the map~$\Upsilon$
    in~\cref{para-upsilon} below respects the tensor product.

    There is an element
    \begin{equation*}
        \bbL^{1/2} = 1 - [\upmu_2]^\muhat \in \Kvarm (\bbK) \ ,
    \end{equation*}
    where $\muhat$ acts on $\upmu_2$ non-trivially.
    It satisfies $(\bbL^{1/2})^2 = \bbL$.
    We also write $\bbL^{-1/2} = \bbL^{-1} \cdot \bbL^{1/2} \in \Mm (\bbK)$.

    There are natural maps
    \begin{equation*}
        \iota^\muhat \colon \Kvar (\calX) \longrightarrow \Kvarm (\calX) \ , \quad
        \iota^\muhat \colon \upM (\calX) \longrightarrow \Mm (\calX) \ , \quad
        \iota^\muhat \colon \Mhat (\calX) \longrightarrow \Mmhat (\calX) \ ,
    \end{equation*}
    given on generators by $[Z] \mapsto [Z]$,
    with the trivial $\muhat$-action on $Z$.
    They are $\Kvar (\bbK)$-, $\upM (\bbK)$-, and $\Mhat (\bbK)$-algebra homomorphisms,
    respectively.

    One can define pullback and pushforward maps
    on $\Kvarm (-)$, $\Mm (-)$, and $\Mmhat (-)$,
    similar to the case of $\Kvar (-)$, $\upM (-)$, and $\Mhat (-)$.
    They satisfy the base change and projection formulae in~\cref{para-mot-bc-pf},
    and the principal bundle relation in~\cref{para-pr-bun-rel}.

    There is also the Euler characteristic map
    $\chi \colon \Mm (\bbK) \to \bbZ$, which is a ring homomorphism,
    defined by taking the Euler characteristic of the underlying non-monodromic motive.
    From this, one can define analogues of operations in
    \cref{para-euler-char,para-cf} for $\Mm (-)$ and $\Mmhat (-)$.
    In particular, we have $\chi (\bbL^{1/2}) = -1$.
\end{para}

\begin{para}[Motives of double covers]
    \label{para-upsilon}
    For a principal $\upmu_2$-bundle $\calP \to \calX$,
    there is a class
    \begin{equation*}
        \Upsilon (\calP) = \bbL^{-1/2} \cdot ([\calX] - [\calP]^\muhat)
        \in \Mmhat (\calX) \ ,
    \end{equation*}
    where $\muhat$ acts on $\calP$ via the $\upmu_2$-action,
    and the class $[\calP]^\muhat$ is defined similarly to
    \cref{para-mot-of-stacks} when $\calP$ is not a variety.
    See \cite[Definitions~5.5 and~5.13]{BBBBJ2015}
    for more details.

    Note that $\Upsilon$~commutes with pullbacks by definition.
    Also, we have the relation
    \begin{equation}
        \Upsilon (\calP_1 \otimes \calP_2) =
        \Upsilon (\calP_1) \cdot \Upsilon (\calP_2)
    \end{equation}
    for principal $\upmu_2$-bundles $\calP_1, \calP_2 \to \calX$,
    where $\calP_1 \otimes \calP_2$ is also a principal $\upmu_2$-bundle.
    See \cite[Definition~5.5 and~5.13]{BBBBJ2015}.
\end{para}

\subsection{Descent of motives}
\label{subsec-mot-desc}

\begin{para}
    We now discuss descent properties of
    the rings of motives defined above.
    While constructible functions $\CF (-)$ descend
    under any reasonable topology,
    descent for rings of motives such as~$\Mhat (-)$
    is more restrictive.
    For example, pulling back along the double cover
    $(-)^2 \colon \Gm \to \Gm$
    is \emph{not} injective on motives,
    since the class of the trivial double cover $\Gm \times \upmu_2 \to \Gm$
    and the non-trivial double cover $\Gm \to \Gm$
    get identified after pulling back.
    Therefore, rings of motives do \emph{not} satisfy \'etale descent.

    However, we show in \cref{thm-nisn-desc} below
    that these rings of motives do satisfy descent
    under the Nisnevich topology.
\end{para}

\begin{para}[The Nisnevich topology]
    \label{para-nisn-top}
    Recall that for an algebraic space~$X$,
    a \emph{Nisnevich cover} of~$X$
    is a family of \'etale morphisms $(f_i \colon X_i \to X)_{i \in I}$,
    such that for each point $x \in X$,
    there exists $i \in I$ and a point $x' \in X_i$,
    such that $f_i (x') = x$,
    and $f_i$ induces an isomorphism on
    residue fields at~$x'$ and~$x$.

    Let $\calX$ be an algebraic stack.
    Define a \emph{Nisnevich cover} of~$\calX$
    to be a representable \'etale cover $(f_i \colon \calX_i \to \calX)_{i \in I}$
    such that its base change to any algebraic space
    is a Nisnevich cover of algebraic spaces.
    See also
    Choudhury--Deshmukh--Hogadi~\cite[Definition~1.2~ff.]{ChoudhuryDeshmukhHogadi2023}.

    For example, for an integer $n > 1$,
    the morphism $* \to [* / \upmu_n]$ is \emph{not} a Nisnevich cover,
    since its base change $\Gm \to \Gm$, $t \mapsto t^n$
    is not a Nisnevich cover.

    Algebraic spaces over~$\bbK$ (assumed locally of finite type)
    admit Nisnevich covers by affine $\bbK$\=/varieties,
    as can be deduced from Knutson~\cite[II, Theorem~6.4]{Knutson1971}.
\end{para}

\begin{theorem}
    \label{thm-nisn-desc}
    Let $\calX$ be an algebraic stack over~$\bbK$
    with affine stabilizers,
    and let $(f_i \colon \calX_i \to \calX)_{i \in I}$ be a Nisnevich cover.
    Then the map
    \begin{equation*}
        (f_i^*)_{i \in I} \colon
        \Mhat (\calX) \longrightarrow
        \operatorname{eq} \biggl(
            \prod_{i \in I} \Mhat (\calX_i)
            \longrightrightarrows
            \prod_{i, j \in I} \Mhat (\calX_i \times_{\calX} \calX_j)
        \biggr)
    \end{equation*}
    is an isomorphism,
    where the right-hand side is the equalizer of
    the two maps induced by pulling back along projections from each
    $\calX_i \times_{\calX} \calX_j$ to $\calX_i$ and $\calX_j$, respectively.

    The same also holds for $\smash{\Mmhat (-)}$
    in place of $\smash{\Mhat (-)}$.
    Moreover, if\/ $\calX$ is an algebraic space,
    then the same holds for
    $\Kvar (-),$ $\upM (-),$ $\Kvarm (-),$ and $\Mm (-)$.
\end{theorem}

\begin{proof}
    We write down the proof for $\Mhat (-)$,
    and the other cases are similar.

    We first consider the case when $\calX$ is an algebraic space.
    In this case, one can stratify $\calX$
    into locally closed pieces $S_k \subset \calX$, such that the map
    $\coprod_i \calX_i \to \calX$ admits a section $s_k$ over each $S_k$.
    After a base change to each $S_k$, we can assume that
    $\coprod_i \calX_i \to \calX$ admits a global section,
    in which case the result is clear.

    For the general case,
    by Kresch~\cite[Proposition~3.5.9]{Kresch1999},
    $\calX$ can be stratified by quotient stacks
    of the form~$[U / G]$,
    where $U$ is a quasi-projective $\bbK$\=/variety
    acted on by $G \simeq \GL (n)$ for some $n$.
    Therefore, we may assume that $\calX = [U / G]$ is of this form.
    Let $\pi \colon U \to [U / G]$ be the projection.
    Then for all $a \in \Mhat ([U / G])$,
    we have $a = [G]^{-1} \cdot \pi_! \circ \pi^* (a)$,
    so that $\pi^*$ is injective.
    Its image consists of elements $\tilde{a} \in \Mhat (U)$
    such that $\pi^* \circ \pi_! (\tilde{a}) = [G] \cdot \tilde{a}$.
    We call such elements \emph{$G$-invariant}.
    In other words, we may identify $\Mhat ([U / G])$
    with the subring of $\Mhat (U)$ consisting of $G$-invariant elements.
    Writing $U_i = U \times_{\calX} \calX_i$, it suffices to show that
    $\Mhat (U) \simto \operatorname{eq} \bigl(
        \prod_{i \in I} \Mhat (U_i) \rightrightarrows
        \prod_{i, j \in I} \Mhat (U_i \times_U U_j) \bigr)$,
    since taking $G$-invariant elements on both sides
    gives the desired result.
    We are now reduced to the already known case of algebraic spaces.
\end{proof}

\begin{para}[Quotient stacks and fundamental stacks]
    \label{para-quot-st}
    We now discuss classes of algebraic stacks that
    can be covered by quotient stacks,
    which will be used in the sequel.
    See \cref{para-quot-st-eg} below for
    examples of stacks satisfying these properties.

    A \emph{quotient stack} over~$\bbK$
    is an algebraic stack of the form $[U / G]$,
    where $U$ is an algebraic space over~$\bbK$ acted on by
    an algebraic group~$G \simeq \GL (n)$ for some~$n$.
    Equivalently, one can allow $G$ to be any linear algebraic group,
    since if one chooses an embedding $G \hookrightarrow \GL (n)$,
    then $[U / G] \simeq [(U \times^G \GL (n)) / \GL (n)]$.

    A \emph{fundamental stack} over~$\bbK$
    is a quotient stack of the form $[U / G]$,
    where $U$ is an affine $\bbK$\=/variety acted on by
    an algebraic group~$G \simeq \GL (n)$ for some~$n$.
    Equivalently, one can allow $G$ to be any reductive group,
    by a similar argument as above;
    see the proof of \cite[Corollary~4.14]{AlperHallRydh2020}.
    This terminology is from Alper--Hall--Rydh~\cite{AlperHallRydh2023}.

    An algebraic stack over~$\bbK$ is \emph{Nisnevich locally a quotient stack},
    if it admits a Nisnevich cover by quotient stacks.
    This class of stacks is also discussed in
    Choudhury--Deshmukh--Hogadi~\cite{ChoudhuryDeshmukhHogadi2023},
    where they are called \emph{cd-quotient stacks}.

    An algebraic stack over~$\bbK$ is \emph{Nisnevich locally fundamental},
    if it admits a Nisnevich cover by fundamental stacks.
    This implies being Nisnevich locally a quotient stack.

    Similarly, we define stacks that are
    \emph{\'etale locally a quotient stack} or
    \emph{\'etale locally fundamental},
    by requiring representable \'etale covers
    with the corresponding properties.

    These properties are satisfied by a large class of stacks,
    which we discuss below.
\end{para}

\begin{para}[Local structure theorems]
    \label{para-quot-st-eg}
    A series of local structure theorems for algebraic stacks
    by Alper--Hall--Rydh~\cite{AlperHallRydh2020,AlperHallRydh2023}
    can be used to produce covers of algebraic stacks
    by quotient and fundamental stacks.
\end{para}

\begin{theorem*}
    Let $\calX$ be an algebraic stack with affine stabilizers.

    \begin{enumerate}[beginpenalty=10000]
        \item
            If every point of\/ $\calX$ specializes to a closed point,
            and closed points of\/ $\calX$ have reductive stabilizers,
            then $\calX$ is \'etale locally fundamental.

        \item
            If\/ $\calX$ admits a good moduli space
            in the sense of Alper~\textnormal{\cite{Alper2013}},
            then $\calX$ is Nisnevich locally fundamental.
    \end{enumerate}
\end{theorem*}

The first result follows from~\cite[Theorem~1.1]{AlperHallRydh2020},
and is stated in Alper--Halpern-Leistner--Heinloth~%
\cite[Remarks~2.6 and~2.7]{AlperHalpernLeistnerHeinloth2023}.
The second result is~\cite[Theorem~6.1]{AlperHallRydh2023}.

\subsection{Motivic vanishing cycles for schemes}
\label{subsec-nearby-sch}

\begin{para}[Idea]
    The motivic nearby and vanishing cycle maps considered below
    are based on the idea of \emph{Milnor fibres},
    which is perhaps more straightforward to explain in the analytic setting.

    For this purpose, let $X$ be a complex manifold,
    with a smooth metric $d \colon X \times X \to \bbR_{\geq 0}$,
    and let $f \colon X \to \bbC$ be a holomorphic function.
    Let $x \in X$ be a point such that $f (x) = 0$.
    Let $0 < \delta \ll \varepsilon \ll 1$ be small positive numbers, and consider the map
    \begin{equation*}
        X^\times_{\delta, \varepsilon} (x) = B_{\varepsilon} (x) \cap f^{-1} (D^\times_\delta)
        \overset{f}{\longrightarrow} D^\times_\delta \ ,
    \end{equation*}
    where $D^\times_\delta = \{ z \in \bbC \mid 0 < |z| < \delta \}$.
    This map is a topological fibration,
    and its fibre $\MF_f (x)$ is called the \emph{Milnor fibre} of~$f$ at~$x$.
    The cohomology of $\MF_f (x)$ carries the action of the
    \emph{monodromy operator} induced by this fibration.

    We will consider below the motivic analogue of the above construction,
    with $X$ a smooth variety over a field~$\bbK$.
    The Milnor fibre is then replaced by the \emph{motivic Milnor fibre}
    of Denef--Loeser~\cite{DenefLoeser1998,DenefLoeser2001,DenefLoeser2002},
    which is a monodromic motive on $X_0$.
    This construction is closely related to Donaldson--Thomas theory;
    see Behrend~\cite{Behrend2009}, Joyce--Song~\cite[\S4]{JoyceSong2012},
    and Kontsevich--Soibelman~\cite[\S4]{KontsevichSoibelman2008}
    for more details.
\end{para}

\begin{para}[The motivic Milnor fibre]
    \label{para-def-mf}
    Let~$X$ be a smooth, irreducible~$\bbK$\=/variety,
    and let $f \colon X \to \smash{\bbA^1}$ be a morphism.
    Write $X_0 = \smash{f^{-1}} (0)$.
    Following Denef--Loeser~\cite{DenefLoeser1998,DenefLoeser2001,DenefLoeser2002},
    we define the \emph{motivic Milnor fibre} of~$f$,
    which is an element
    \begin{equation*}
        \MF_f \in \Mm (X_0) \ ,
    \end{equation*}
    as follows.

    If $f$ is constant, define $\smash{\MF_f} = 0$.
    Otherwise, choose a resolution $\pi \colon \tilde{X} \to X$ of~$f$,
    in the sense that $\tilde{X}$ is a smooth, irreducible $\bbK$\=/variety,
    $\pi$ is a proper morphism that restricts to
    an isomorphism on~$\pi^{-1} (X \setminus X_0)$,
    and~$\pi^{-1} (X_0)$ is a simple normal crossings divisor in~$\tilde{X}$.
    See, for example, Koll\'ar~\cite{Kollar2007}
    for the existence and properties of such resolutions.

    Let $\smash{(E_i)_{i \in J}}$ be the irreducible components of~$\pi^{-1} (X_0)$,
    and write $N_i$ for the multiplicity of~$E_i$
    in the divisor of $f \circ \pi$ on~$\tilde{X}$.
    For a non-empty subset $I \subset J$, write
    $E^\circ_I = \bigcap_{i \in I} E_i \setminus \bigcup_{i \notin I} E_i$.
    Let $m_I = \gcd_{i \in I} N_i$,
    and define an $m_I$-fold cover $\tilde{E}^\circ_I \to E^\circ_I$ as follows.
    For each open set $U \subset \tilde{X}$ such that
    $f \circ \pi = u v^{m_I}$ on~$U$ for
    $u \colon U \to \bbA^1 \setminus \{ 0 \}$ and $v \colon U \to \bbA^1$,
    define the restriction of~$\tilde{E}^\circ_I$
    on~$E^\circ_I \cap U$ as
    \begin{equation}
        \label{eq-def-e-tilde-circ}
        \tilde{E}^\circ_I |_{E^\circ_I \cap U} = \bigl\{
            (z, y) \in \bbA^1 \times (E^\circ_I \cap U)
            \bigm| z^{m_I} = u^{-1}
        \bigr\} \ .
    \end{equation}
    Since $E^\circ_I$ can be covered by such open sets~$U$,
    \cref{eq-def-e-tilde-circ} can be glued together
    to obtain a cover $\tilde{E}^\circ_I \to E^\circ_I$,
    with a natural $\upmu_{\smash{m_I}}$-action given by scaling the $z$-coordinate,
    which induces a $\muhat$-action on $\tilde{E}^\circ_I$.
    The motivic Milnor fibre $\smash{\MF_f}$ is then given by
    \begin{equation}
        \label{eq-def-mf}
        \MF_f = \sum_{\varnothing \neq I \subset J} {} (1 - \bbL)^{|I| - 1} \,
        [\tilde{E}^\circ_I]^\muhat \ .
    \end{equation}
    It can be shown \cite[Definition~3.8]{DenefLoeser2001}
    that this is independent of the choice of the resolution~$\pi$.
\end{para}

\begin{para}[Nearby and vanishing cycles]
    \label{para-def-nearby}
    Let $X$ be a $\bbK$\=/variety,
    and let $f \colon X \to \smash{\bbA^1}$ be a morphism.
    Write $X_0 = \smash{f^{-1}} (0)$.

    Define the \emph{nearby cycle map} of~$f$, denoted by
    \begin{equation*}
        \Psi_f \colon \upM (X) \longrightarrow \Mm (X_0) \ ,
    \end{equation*}
    to be the unique $\upM (\bbK)$-linear map such that
    for any smooth, irreducible $\bbK$\=/variety~$Z$
    and any proper morphism $\rho \colon Z \to X$, we have
    \begin{equation*}
        \Psi_f ([Z]) = (\rho_0)_! \, (\MF_{f \circ \rho}) \in \Mm (X_0) \ ,
    \end{equation*}
    where $\rho_0 \colon Z_0 \to X_0$ is the restriction of $\rho$ to
    $Z_0 = (f \circ \rho)^{-1} (0)$,
    and $\MF_{f \circ \rho} \in \Mm (Z_0)$
    is the motivic Milnor fibre of $f \circ \rho$.
    It follows from Bittner~\cite[Claim~8.2]{Bittner2005}
    that the map~$\Psi_f$ is well-defined.

    Define the \emph{vanishing cycle map} of~$f$ to be the map
    \begin{equation*}
        \Phi_f = \Psi_f - \iota^\muhat \circ i^* \colon
        \upM (X) \longrightarrow \Mm (X_0) \ ,
    \end{equation*}
    where $i \colon X_0 \hookrightarrow X$ is the inclusion,
    and $\iota^\muhat \colon \upM (X_0) \hookrightarrow \Mm (X_0)$
    is as in~\cref{para-mon-mot}.
\end{para}

\begin{para}[For algebraic spaces]
    \label{para-nearby-alg-sp}
    We now generalize the above construction
    from varieties to algebraic spaces.

    As in Bittner~\cite[Theorem~8.4]{Bittner2005},
    the nearby and vanishing cycle maps are compatible
    with pulling back along smooth morphisms.
    In particular, they define morphisms
    $\Psi, \Phi \colon \upM (-) \to \Mm ((-)_0)$
    of sheaves on the category of $\bbK$\=/varieties over~$\bbA^1$,
    with the Nisnevich topology.
    Since algebraic spaces admit Nisnevich covers
    by affine $\bbK$\=/varieties, as mentioned in \cref{para-nisn-top},
    these morphisms of sheaves induce maps
    on their evaluations on algebraic spaces over~$\bbK$.

    In other words, for an algebraic space~$X$ over~$\bbK$
    and a morphism $f \colon X \to \bbA^1$,
    we have defined nearby and vanishing cycle maps
    \begin{equation*}
        \Psi_f, \Phi_f \colon \upM (X) \longrightarrow \Mm (X_0) \ .
    \end{equation*}
    We state some of their properties below.
\end{para}

\begin{theorem}
    \label{thm-sch-nearby-properties}
    Let $X, Y$ be algebraic spaces over~$\bbK$.

    \begin{enumerate}
        \item \label{itm-sch-proper-pushforward}
            Let $g \colon Y \to X$ be a proper morphism,
            and $f \colon X \to \bbA^1$ a morphism.
            Then we have a commutative diagram
            \begin{equation}
                \begin{tikzcd}
                    \upM (Y) \ar[r, "g_!"] \ar[d, "\Psi_{f \circ g}"'] &
                    \upM (X) \ar[d, "\Psi_f"]
                    \\
                    \Mm (Y_0) \ar[r, "g_!"] &
                    \Mm (X_0) \rlap{\textnormal{ .}}
                \end{tikzcd}
            \end{equation}

        \item \label{itm-sch-smooth-pullback}
            Let $g \colon Y \to X$ be a smooth morphism,
            and $f \colon X \to \bbA^1$ a morphism.
            Then we have a commutative diagram
            \begin{equation}
                \begin{tikzcd}
                    \upM (X) \ar[r, "g^*"] \ar[d, "\Psi_f"'] &
                    \upM (Y) \ar[d, "\Psi_{f \circ g}"]
                    \\
                    \Mm (X_0) \ar[r, "g^*"] &
                    \Mm (Y_0) \rlap{\textnormal{ .}}
                \end{tikzcd}
            \end{equation}
    \end{enumerate}
\end{theorem}

\begin{proof}
    The case when $X$ and $Y$ are $\bbK$\=/varieties
    was proved in Bittner~\cite[Theorem~8.4]{Bittner2005}.
    The verification of \cref{itm-sch-smooth-pullback}
    for algebraic spaces is completely formal,
    by passing to Nisnevich covers by $\bbK$\=/varieties.

    We now prove \cref{itm-sch-proper-pushforward}
    for algebraic spaces.
    Again, passing to a Nisnevich cover,
    we may assume that $X$ is a $\bbK$\=/variety.
    We claim that $\Kvar (Y)$ is spanned over~$\bbZ$
    by classes $[Z]$ of proper morphisms $Z \to Y$,
    where $Z$ is a smooth $\bbK$\=/variety.
    Indeed, let $u \colon U \to Y$ be an arbitrary morphism,
    where $U$ is an integral $\bbK$\=/variety.
    By Nagata compactification,
    as in Conrad--Lieblich--Olsson~%
    \cite[Theorem~1.2.1]{ConradLieblichOlsson2012},
    $u$ can be factored as a dense open immersion $U \hookrightarrow V$
    followed by a proper morphism $V \to Y$,
    where $V$ is an integral algebraic space over~$\bbK$.
    By Chow's lemma for algebraic spaces,
    as in Knutson~\cite[IV, Theorem~3.1]{Knutson1971},
    there exists a $\bbK$\=/variety $W$
    and a projective birational morphism $W \to V$.
    Applying a resolution of singularities,
    we may assume that $W$ is smooth.
    Now $W \to Y$ is proper,
    and the difference $[W] - [U]$
    is a sum of lower dimensional classes.
    An induction on the dimension of~$U$ verifies the claim.

    Now, let $h \colon Z \to Y$ be a proper morphism,
    where $Z$ is a smooth $\bbK$\=/variety.
    Passing to a Nisnevich cover of~$Y$ by $\bbK$\=/varieties,
    one can show that
    $\Psi_{f \circ g} ([Z]) = h_! (\MF_{f \circ g \circ h})$.
    On the other hand, we have
    $\Psi_f ([Z]) = (g \circ h)_! (\MF_{f \circ g \circ h})$ by definition.
    This completes the proof since
    such classes $[Z]$ span~$\Kvar (Y)$,
    so they also span~$\upM (Y)$ over~$\upM (\bbK)$.
\end{proof}

\subsection{Motivic vanishing cycles for stacks}
\label{subsec-nearby-stacks}

\begin{para}[Assumptions on the stack]
    \label{para-nearby-assumptions}
    From now on,
    we assume that $\calX$ is an algebraic stack over~$\bbK$ that is
    Nisnevich locally a quotient stack
    in the sense of \cref{para-quot-st}.

    Note that this assumption is weaker than that in
    Ben-Bassat--Brav--Bussi--Joyce~\cite[\S5]{BBBBJ2015},
    where the stack was assumed to be
    Zariski locally a quotient stack.
\end{para}

\begin{theorem}
    \label{thm-nearby-stacks}
    Let $\calX$ be as in~\cref{para-nearby-assumptions},
    and let $f \colon \calX \to \bbA^1$ be a morphism.
    Write $\calX_0 = \smash{f^{-1}} (0)$.
    Then there is a unique $\Mhat (\bbK)$-linear map
    \begin{equation*}
        \Psi_f \colon \Mhat (\calX) \longrightarrow \Mmhat (\calX_0) \ ,
    \end{equation*}
    called the \emph{nearby cycle map} of\/~$f$,
    such that for any $\bbK$\=/variety~$Y$
    and any smooth morphism $g \colon Y \to \calX$,
    we have a commutative diagram
    \begin{equation}
        \label{eq-nearby-stacks-cond}
        \begin{tikzcd}
            \Mhat (\calX) \ar[d, "\Psi_f"'] \ar[r, "g^*"] &
            \Mhat (Y) \ar[d, "\Psi_{f \circ g}"]
            \\
            \Mmhat (\calX_0) \ar[r, "g^*"] &
            \Mmhat (Y_0) \rlap{ ,}
        \end{tikzcd}
    \end{equation}
    where the right-hand map is defined in~\cref{para-def-nearby}.

    We then define the \emph{vanishing cycle map} of\/~$f$ to be the map
    \begin{equation*}
        \Phi_f = \Psi_f - \iota^\muhat \circ i^* \colon
        \Mhat (\calX) \longrightarrow \Mmhat (\calX_0) \ ,
    \end{equation*}
    where $i \colon \calX_0 \hookrightarrow \calX$
    and $\iota^\muhat \colon \Mhat (\calX_0) \hookrightarrow \Mmhat (\calX_0)$
    are the inclusions.
\end{theorem}

\begin{proof}
    \allowdisplaybreaks
    Let $(j_i \colon \calX_i \to \calX)_{i \in I}$ be a Nisnevich cover,
    where each $\calX_i \simeq [U_i / G_i]$,
    with $U_i$ an algebraic space over~$\bbK$, acted on by a group
    $G_i \simeq \GL (n_i)$ for some~$n_i$.
    Let $\pi_i \colon U_i \to \calX_i$ be the projection.

    First, note that the condition on~$\smash{\Psi_f}$
    implies that the same condition holds when~$Y$ is an algebraic space,
    with the right-hand map in~\cref{eq-nearby-stacks-cond}
    defined in~\cref{para-nearby-alg-sp}.
    This can be seen by passing to
    a Nisnevich cover of~$Y$ by $\bbK$\=/varieties,
    and applying \cref{thm-nisn-desc} to this cover.

    To define the map~$\smash{\Psi_f}$, by \cref{thm-nisn-desc},
    it is enough to define it on each~$\calX_i$,
    and then verify that they agree on overlaps.
    Let $a \in \Mhat (\calX)$ be an element.
    We define the element $\smash{\Psi_f (a)} \in \Mmhat (\calX_0)$
    by giving its pullbacks
    $\smash{\Psi_f (a)_i} = j_i^* \circ \smash{\Psi_f (a)}
    \in \Mmhat (\calX_{i, 0})$ for each~$i$,
    where $\calX_{i, 0} = \calX_i \times_{\calX} \calX_0$.
    The condition on $\smash{\Psi_f}$ forces
    \begin{align*}
        \Psi_f (a)_i =
        j_i^* \circ \Psi_f (a)
        & =
        [G_i]^{-1} \cdot (\pi_i)_! \circ \pi_i^* \circ j_i^* \circ \Psi_f (a)
        \\* & =
        [G_i]^{-1} \cdot (\pi_i)_! \circ \Psi_{f \circ j_i \circ \pi_i} \circ \pi_i^* \circ j_i^* (a) \ ,
        \numberthis
    \end{align*}
    where $[G_i] \in \Mhat (\bbK)$ is the class of~$G_i$
    and is invertible in $\Mhat (\bbK)$,
    and we applied~\cref{eq-pr-bun-rel} to $\pi_i$,
    using the fact that~$G_i$ is special.
    This shows that if the map $\smash{\Psi_f}$ exists,
    then it is unique.

    To check that the elements~$\Psi_f (a)_i$ agree on overlaps,
    let $1, 2 \in I$ be two indices, and form the pullback squares
    \begin{equation}
        \begin{tikzcd}[column sep={5em,between origins}]
            U''
            \ar[r, "\pi''_2"] \ar[d, "\pi''_1"']
            \ar[dr, phantom, pos=.3, start anchor=center, end anchor=center, "\ulcorner"] &
            U'_1
            \ar[r, "j''_2"] \ar[d, "\pi'_1"']
            \ar[dr, phantom, pos=.3, start anchor=center, end anchor=center, "\ulcorner"] &
            U_1 \ar[d, "\pi_1"]
            \\
            U'_2
            \ar[r, "\pi'_2"] \ar[d, "j''_1"']
            \ar[dr, phantom, pos=.3, start anchor=center, end anchor=center, "\ulcorner"] &
            \calX_{1, 2}
            \ar[r, "j'_2"] \ar[d, "j'_1"']
            \ar[dr, phantom, pos=.3, start anchor=center, end anchor=center, "\ulcorner"] &
            \calX_1 \ar[d, "j_1"]
            \\
            U_2 \ar[r, "\pi_2"] &
            \calX_2 \ar[r, "j_2"] &
            \calX \rlap{ ,}
        \end{tikzcd}
    \end{equation}
    where $U'_1$, $U'_2$, $U''$ are algebraic spaces.
    We need to show that
    $(j'_2)^* (\smash{\Psi_f (a)_1}) = (j'_1)^* (\smash{\Psi_f (a)_2})$.
    We have
    \begin{align*}
        & \phantom{{} = {}}
        (j'_2)^* (\Psi_f (a)_1)
        \\* & =
        [G_1]^{-1} \cdot
        (j'_2)^* \circ (\pi_1)_! \circ
        \Psi_{f \circ j_1 \circ \pi_1} \circ (j_1 \circ \pi_1)^* (a)
        \\ & =
        [G_1]^{-1} \cdot
        (\pi'_1)_! \circ (j''_2)^* \circ
        \Psi_{f \circ j_1 \circ \pi_1} \circ (j_1 \circ \pi_1)^* (a)
        \\ & =
        [G_1]^{-1} \cdot
        (\pi'_1)_! \circ \Psi_{f \circ j_1 \circ \pi_1 \circ j''_2} \circ
        (j_1 \circ \pi_1 \circ j''_2)^* (a)
        \\ & =
        [G_1]^{-1} \cdot [G_2]^{-1} \cdot
        (\pi'_1)_! \circ (\pi''_2)_! \circ (\pi''_2)^* \circ
        \Psi_{f \circ j_1 \circ \pi_1 \circ j''_2} \circ
        (j_1 \circ \pi_1 \circ j''_2)^* (a)
        \\* & =
        [G_1]^{-1} \cdot [G_2]^{-1} \cdot
        (\pi'_1 \circ \pi''_2)_! \circ
        \Psi_{f \circ j_1 \circ \pi_1 \circ j''_2 \circ \pi''_2} \circ
        (j_1 \circ \pi_1 \circ j''_2 \circ \pi''_2)^* (a) \ ,
        \numberthis
    \end{align*}
    where we applied \cref{eq-mot-bc} in the second step,
    \cref{thm-sch-nearby-properties}~\ref{itm-sch-smooth-pullback}
    in the third and fifth steps,
    and \cref{para-pr-bun-rel} in the fourth step.
    This expression is now symmetric in the indices $1$ and~$2$,
    so the element $\smash{\Psi_f} (a)$ is well-defined.

    It now remains to show that the map~$\smash{\Psi_f}$
    satisfies the required condition.
    Let $Y$ be a $\bbK$\=/variety and
    $\pi \colon Y \to \calX$ a smooth morphism.
    For each $i \in I$, write $Y_i = Y \times_{\calX} \calX_i$.
    Then $(k_i \colon Y_i \to Y)_{i \in I}$
    is a Nisnevich cover by algebraic spaces,
    and it suffices to show that
    \begin{equation}
        \label{eq-pf-nearby-wd}
        k_i^* \circ g^* \circ \Psi_f = k_i^* \circ \Psi_{f \circ g} \circ g^*
    \end{equation}
    for each~$i$. Consider the diagram
    \begin{equation}
        \begin{tikzcd}
            V_i \ar[r, "\rho_i"] \ar[d, "g'_i"']
            \ar[dr, phantom, pos=.2, "\ulcorner"] &
            Y_i \ar[r, "k_i"] \ar[d, "g_i"']
            \ar[dr, phantom, pos=.2, "\ulcorner"] &
            Y \ar[d, "g"] \\
            U_i \ar[r, "\pi_i"] &
            \calX_i \ar[r, "j_i"] &
            \calX \rlap{ ,}
        \end{tikzcd}
    \end{equation}
    where all squares are pullback squares.
    In particular, $\rho_i$ is a principal $G_i$-bundle.
    For any $a \in \Mhat (\calX)$, we have
    \begin{align*}
        & \phantom{{} = {}}
        k_i^* \circ g^* \circ \Psi_f (a)
        \\* & =
        g_i^* (\Psi_f (a)_i)
        \\ & =
        [G_i]^{-1} \cdot
        g_i^* \circ (\pi_i)_! \circ \Psi_{f \circ j_i \circ \pi_i} \circ (j_i \circ \pi_i)^* (a)
        \\ & =
        [G_i]^{-1} \cdot
        (\rho_i)_! \circ (g'_i)^* \circ \Psi_{f \circ j_i \circ \pi_i} \circ (j_i \circ \pi_i)^* (a)
        \\ & =
        [G_i]^{-1} \cdot
        (\rho_i)_! \circ \Psi_{f \circ j_i \circ \pi_i \circ g'_i} \circ (j_i \circ \pi_i \circ g'_i)^* (a)
        \\ & =
        [G_i]^{-1} \cdot
        (\rho_i)_! \circ \Psi_{f \circ g \circ k_i \circ \rho_i} \circ (g \circ k_i \circ \rho_i)^* (a)
        \\ & =
        [G_i]^{-1} \cdot
        (\rho_i)_! \circ \rho_i^* \circ \Psi_{f \circ g \circ k_i} \circ (g \circ k_i)^* (a)
        \\ & =
        \Psi_{f \circ g \circ k_i} \circ (g \circ k_i)^* (a)
        \\* & =
        k_i^* \circ \Psi_{f \circ g} \circ g^* (a) \ ,
        \numberthis
    \end{align*}
    where we applied \cref{eq-mot-bc} in the third step,
    \cref{thm-sch-nearby-properties}~\ref{itm-sch-smooth-pullback}
    in the fourth, sixth, and eighth steps,
    and \cref{para-pr-bun-rel} in the seventh step.
    This proves the desired identity~\cref{eq-pf-nearby-wd}.
\end{proof}

\begin{para}[The motivic Milnor fibre]
    \label{para-milnor-stacks}
    Let $\calX$ be as in~\cref{para-nearby-assumptions},
    and let $f \colon \calX \to \smash{\bbA^1}$ be a morphism.
    Write $\calX_0 = \smash{f^{-1}} (0)$.
    The \emph{motivic Milnor fibre} of~$f$ is the element
    \begin{equation*}
        \MF_f = \Psi_f ([\calX]) \in \Mmhat (\calX_0) \ .
    \end{equation*}
    When $\calX$ is smooth,
    this is closely related to the construction of
    Ben-Bassat--Brav--Bussi--Joyce~\cite[\S5.4]{BBBBJ2015},
    which we will further discuss and generalize
    in \cref{para-nu-mot-st} below.
    The main difference is that
    the latter construction starts from
    the critical locus of~$f$ instead of~$\calX$,
    and can be generalized to stacks glued from such critical loci;
    it uses~$\Phi_{\smash{f}}$ instead of~$\Psi_{\smash{f}}$,
    and introduces a twist by a power of~$\bbL^{1/2}$ to make gluing possible.

    We relate this to the description of the
    motivic Milnor fibre for schemes in~\cref{para-def-mf}.
    Suppose that we are given a \emph{resolution} of~$f$,
    which is a representable proper morphism $\pi \colon \tilde{\calX} \to \calX$,
    such that it restricts to an isomorphism on~$\pi^{-1} (\calX \setminus \calX_0)$,
    and $\pi^{-1} (\calX_0)$ is a simple normal crossings divisor in~$\tilde{\calX}$,
    in the sense that it is so after pulling back along
    smooth morphisms from schemes.
    Let $(\calE_i)_{\smash{i \in J}}$ be the family of
    irreducible components of~$\pi^{-1} (\calX_0)$,
    and define $\calE^\circ_I$ and $\tilde{\calE}^\circ_I$ for non-empty $I \subset J$
    similarly to~\cref{para-def-mf},
    where $\tilde{\calE}^\circ_I$ carries a natural $\muhat$-action.
    We then claim that
    \begin{equation}
        \MF_f = \sum_{\varnothing \neq I \subset J} {} (1 - \bbL)^{|I| - 1} \,
        [\tilde{\calE}^\circ_I]^\muhat \ .
    \end{equation}
    Indeed, this can be shown by a similar argument
    as in the proof of~\cref{thm-nearby-stacks},
    by first passing to a Nisnevich cover by quotient stacks,
    then using the relation~\cref{eq-pr-bun-rel}
    to further reduce to the case of algebraic spaces,
    and finally passing to a Nisnevich cover again
    to reduce to the case of affine varieties.
\end{para}

\begin{theorem}
    \label{thm-nearby-properties}
    Let $\calX, \calY$ be algebraic stacks as in~\cref{para-nearby-assumptions}.

    \begin{enumerate}
        \item \label{itm-proper-pushforward}
            Let $g \colon \calY \to \calX$ be a proper morphism,
            and $f \colon \calX \to \bbA^1$ a morphism.
            Then we have a commutative diagram
            \begin{equation}
                \begin{tikzcd}
                    \Mhat (\calY) \ar[r, "g_!"] \ar[d, "\Psi_{f \circ g}"'] &
                    \Mhat (\calX) \ar[d, "\Psi_f"]
                    \\
                    \Mmhat (\calY_0) \ar[r, "g_!"] &
                    \Mmhat (\calX_0) \rlap{\textnormal{ .}}
                \end{tikzcd}
            \end{equation}

        \item \label{itm-smooth-pullback}
            Let $g \colon \calY \to \calX$ be a smooth morphism,
            and $f \colon \calX \to \bbA^1$ a morphism.
            Then we have a commutative diagram
            \begin{equation}
                \begin{tikzcd}
                    \Mhat (\calX) \ar[r, "g^*"] \ar[d, "\Psi_f"'] &
                    \Mhat (\calY) \ar[d, "\Psi_{f \circ g}"]
                    \\
                    \Mmhat (\calX_0) \ar[r, "g^*"] &
                    \Mmhat (\calY_0) \rlap{\textnormal{ .}}
                \end{tikzcd}
            \end{equation}
            In particular, we have
            $\MF_{f \circ g} = g^* (\MF_f)$.
    \end{enumerate}
\end{theorem}

\begin{proof}
    \allowdisplaybreaks
    For \cref{itm-proper-pushforward},
    we first restrict to the case when $g$ is representable.
    By \cref{thm-nearby-stacks},
    the map $\Psi_f$ is determined by pullbacks along smooth morphisms
    from $\bbK$\=/varieties to~$\calX$,
    so we may assume that $\calX$ is a $\bbK$\=/variety,
    and $\calY$ is an algebraic space that is proper over $\calX$.
    This case is covered by
    \cref{thm-sch-nearby-properties}~\cref{itm-sch-proper-pushforward}.

    For the general case, similarly,
    we may assume that $\calX = X$ is a $\bbK$\=/variety.
    It suffices to show that
    $g_! \circ \Psi_{f \circ g} ([Z]) = \Psi_f \circ g_! ([Z])$
    for smooth $\bbK$\=/varieties $Z$ mapping to~$\calY$,
    as these classes span~$\smash{\Mhat (\calY)}$ over~$\smash{\Mhat (\bbK)}$.
    Since $\calY$ is proper over~$X$ and has affine stabilizers,
    it has finite inertia,
    and admits a coarse space $\pi_{\calY} \colon \calY \to \oline{Y}$
    by the Keel--Mori theorem~\cite{KeelMori1997,Conrad2005}.
    The morphism $\pi_{\calY}$ is a proper universal homeomorphism.

    By Rydh's compactification theorem for representable morphisms
    of Deligne--Mumford stacks~\cite[Theorem~B]{Rydh2011},
    we may choose a relative compactification~$\calZ$ of~$Z$ over~$\calY$,
    such that there is a dense open immersion $i \colon Z \hookrightarrow \calZ$
    and a proper representable morphism $h \colon \calZ \to \calY$.
    In particular, $\calZ$ also has finite inertia,
    and admits a coarse space $\pi_{\calZ} \colon \calZ \to \oline{Z}$,
    which can be seen as a relative compactification of~$Z$ over~$\oline{Y}$.
    We have a commutative diagram
    \begin{equation}
        \begin{tikzcd}[column sep=1em]
            Z \ar[rr, hook, "i"] &&
            \calZ \ar[d, "h"'] \ar[rr, "\pi_{\calZ}"] &&
            \oline{Z} \ar[d, "\oline{h}"]
            \\
            && \calY \ar[rr, "\pi_{\calY}"] \ar[dr, "g"'] &&
            \oline{Y} \vphantom{\overline{0}} \ar[dl, "\overline{g}"]
            \\
            &&& X \rlap{\textnormal{ ,}}
        \end{tikzcd}
    \end{equation}
    where $\overline{g}$ and $\oline{h}$ are the induced morphisms,
    and all morphisms except~$i$ are proper.
    It is then enough to show that
    \begin{equation}
        \label{eq-pf-nearby-proper-pf}
        (\pi_{\calZ})_! \circ \Psi_{f \circ g \circ h} ([Z]) =
        \Psi_{f \circ \overline{g} \circ \oline{h}} \circ (\pi_{\calZ})_! ([Z]) \ ,
    \end{equation}
    since the compatibility with pushing forward along
    $h$ and $\overline{g} \circ \oline{h}$ is covered by the previous case.

    We now apply Bergh--Rydh's \emph{divisorialification theorem}
    \cite[Theorem~A]{BerghRydh2019} to
    a desingularization of the pair~$(\calZ, \calZ \setminus Z)$
    (see, for example, \cite{EncinasVillamayor1998}),
    which gives a representable proper morphism
    $\oline{\calZ} \to \calZ$ that is an isomorphism on the preimage of~$Z$,
    such that $\oline{\calZ} \setminus Z = \calD$
    is a simple normal crossings divisor on~$\oline{\calZ}$,
    with smooth irreducible components $\calD_i \subset \oline{\calZ}$,
    and for each $x \in \oline{\calZ}$,
    writing $I_x = \{ i \in I \mid x \in \calD_i \}$,
    \'etale locally around~$x$, one has
    $\oline{\calZ} \sim \prod_{\smash{i \in I_x}} {} [\bbA^1 / \smash{\upmu_{n_i}}] \times \bbA^{d - |I_x|}$,
    where $d = \dim \oline{\calZ}$,
    each $\smash{\upmu_{n_i}}$ acts on $\bbA^1$ by scaling,
    and $\calD_i$ corresponds to the locus where the $i$-th factor is zero;
    the number $n_i$ is the order of the generic stabilizer of~$\calD_i$.

    From now on, we assume that $\calZ = \oline{\calZ}$,
    since again, pushing forward along the representable morphism
    $\oline{\calZ} \to \calZ$ and the corresponding morphism
    of coarse spaces is already dealt with.

    Now, choose a resolution $\pi \colon \tilde{\calZ} \to \calZ$
    for the morphism $\calZ \to \bbA^1$,
    which is a composition of blow-ups along smooth centres.
    Then $\tilde{\calZ}$ still has the same local description as before.
    The local description implies that the coarse space of $\tilde{\calZ}$,
    denoted~$\tilde{Z}$, is a smooth algebraic space,
    and can be seen as a resolution for the morphism $\oline{Z} \to \bbA^1$.

    For each $i \in I$, let $\tilde{\calD}_i \subset \tilde{\calZ}$
    be the strict transform of~$\calD_i$,
    which is a smooth divisor,
    and let $(\smash{\calE_j} \subset \tilde{\calZ})_{\smash{j \in J}}$ be the family of
    irreducible components of~$\tilde{\calZ}_0$.
    Then by construction, all the divisors $\tilde{\calD}_i, \smash{\calE_j} \subset \tilde{\calZ}$
    have simple normal crossings,
    and $\tilde{\calZ} \setminus \bigcup_{i \in I} \tilde{\calD}_i$ is an algebraic space.
    Let $\tilde{D}_i, \smash{E_j} \subset \tilde{Z}$
    be the corresponding divisors in the coarse spaces.
    For $I' \subset I$, write $\calD_{I'} = \bigcap_{i \in I'} \calD_i$
    and $\tilde{\calD}_{I'} = \bigcap_{i \in I'} \tilde{\calD}_i$,
    with the convention that $\calD_{\varnothing} = \calZ$
    and $\tilde{\calD}_{\varnothing} = \tilde{\calZ}$.
    Then, each $\tilde{\calD}_{I'}$ can be seen as a resolution
    for the morphism $\calD_{I'} \to \bbA^1$.
    By~\cref{para-milnor-stacks}, we have
    \begin{align*}
        (\pi_{\calZ})_! \circ \Psi_{f \circ g \circ h} ([Z]) & =
        \sum_{I' \subset I} {} (-1)^{|I'|} \cdot
        (\pi_{\calZ})_! \circ \Psi_{f \circ g \circ h} ([\calD_{I'}])
        \\* & =
        \sum_{I' \subset I} {} (-1)^{|I'|} \cdot
        \sum_{\varnothing \neq J' \subset J} {} (1 - \bbL)^{|J'| - 1} \,
        [\tilde{\calE}^\circ_{\smash{J'}} \cap \tilde{\calD}_{I'}]^\muhat
        \\ & =
        \sum_{\varnothing \neq J' \subset J} {} (1 - \bbL)^{|J'| - 1} \,
        \Bigl[ \tilde{\calE}^\circ_{\smash{J'}} \Bigm\backslash
        \bigcup_{i \in I} \tilde{\calD}_i \Bigr]^\muhat
        \\ & =
        \sum_{\varnothing \neq J' \subset J} {} (1 - \bbL)^{|J'| - 1} \,
        \Bigl[ \tilde{E}^\circ_{\smash{J'}} \Bigm\backslash
        \bigcup_{i \in I} \tilde{D}_i \Bigr]^\muhat
        \\* & =
        \Psi_{f \circ \overline{g} \circ \oline{h}} \circ (\pi_{\calZ})_! ([Z]) \ ,
    \end{align*}
    where the second last step used the fact that each
    $\tilde{\calE}^\circ_{\smash{J'}} \setminus \bigcup_{i \in I} \tilde{\calD}_i$
    is an algebraic space.

    For \cref{itm-smooth-pullback},
    similarly, the case when $g$ is representable
    follows from \cref{thm-sch-nearby-properties}~\cref{itm-sch-smooth-pullback}.
    For the general case,
    we may assume that $\calX$ is a $\bbK$\=/variety.
    Since $\Psi_{f \circ g}$ is determined by pullbacks
    along smooth morphisms from schemes to~$\calY$,
    we can also assume that $\calY$ is a $\bbK$\=/variety,
    and the result follows from
    \cref{thm-sch-nearby-properties}~\cref{itm-sch-smooth-pullback}.

    The final statement follows from
    applying~\cref{itm-smooth-pullback}
    to the element $[\calX] \in \Mhat (\calX)$.
\end{proof}

\begin{para}[Remark]
    The length of the proof of \cref{thm-nearby-properties}
    is primarily due to the case of pushing forward
    along proper morphisms that are not necessarily representable.
    We will indeed need this general case in the proof of
    one of our main results, \cref{thm-behrend-local},
    where~$g$ will be taken to be a \emph{weighted blow-up}
    in the sense of~\cref{para-wbl}.
\end{para}

\subsection{The motivic Behrend function}
\label{subsec-nu-mot}

\begin{para}[Shifted symplectic and d\=/critical stacks]
    \label{para-d-crit}
    Let $\frX$ be a $(-1)$\=/shifted symplectic stack over~$\bbK$
    (see \cref{para-conventions}),
    and let $\calX = \frX_\mathrm{cl}$ be its classical truncation.
    Assume that $\calX$ is an algebraic stack over~$\bbK$.

    Ben-Bassat--Brav--Bussi--Joyce~\cite[\S3.3]{BBBBJ2015}
    define a \emph{d\=/critical structure} on~$\calX$
    induced from the shifted symplectic structure on~$\frX$,
    so that $\calX$ is a \emph{d\=/critical stack}.
    See there and Joyce~\cite{Joyce2015}
    for the precise definitions.
    For our purposes, it suffices to know the following properties:

    \begin{enumerate}
        \item \label{itm-d-crit-1}
            For a smooth $\bbK$\=/variety~$U$ and a morphism $f \colon U \to \bbA^1$,
            the critical locus $\Crit (f) \subset U$
            carries a canonical d\=/critical structure.
        \item
            d\=/critical structures can be pulled back along
            smooth morphisms of algebraic stacks over~$\bbK$.
        \item
            If a $\bbK$\=/scheme~$X$ carries a d\=/critical structure,
            then it can be covered by open subschemes
            called \emph{critical charts},
            each of which with the induced d\=/critical structure
            has the form $\Crit (f)$ as in~\cref{itm-d-crit-1}.
            We denote such a critical chart by
            $i \colon \Crit (f) \hookrightarrow X$.
    \end{enumerate}
\end{para}

\begin{para}[Orientations]
    \label{para-ori}
    Let $\frX$ be an $s$\=/shifted symplectic stack over~$\bbK$,
    where $s$~is odd.
    An \emph{orientation} of~$\frX$
    is a line bundle $K_{\frX}^{\smash{1/2}} \to \frX$,
    together with an isomorphism
    $(K_{\frX}^{\smash{1/2}})^{\smash{\otimes 2}} \simeq K_{\frX}$,
    where $K_{\frX}$ is the \emph{canonical bundle} of~$\frX$,
    defined as the determinant line bundle of the cotangent complex of~$\frX$.

    Note the unfortunate clash of terminology with the unrelated notion of
    \emph{orientations} in the sense of Pantev--To\"en--Vaqui\'e--Vezzosi~%
    \cite[Definition~2.4]{PTVV2013}.
    The latter notion will not be used in this paper
    except in the proof of \cref{thm-filt-lag}.

    Now let $s = -1$,
    and let $\calX$ be the associated d\=/critical stack of~$\frX$,
    as in \cref{para-d-crit}.
    Ben-Bassat--Brav--Bussi--Joyce~\cite[Theorem~3.18]{BBBBJ2015}
    show that the restriction~$K_{\frX} |_{\calX^\red}$
    is determined by the d\=/critical structure on~$\calX$,
    where $\calX^\red$ is the reduction of~$\calX$.
    We denote this restriction simply by~$K_{\calX}$,
    and call it the \emph{canonical bundle} of the d\=/critical stack~$\calX$.
    Similarly, as in \cite[Definition~3.6]{BBBBJ2015},
    an \emph{orientation} of a d\=/critical stack~$\calX$
    is a line bundle $K_{\calX}^{\smash{1/2}} \to \calX^\red$,
    together with an isomorphism
    $(K_{\calX}^{\smash{1/2}})^{\smash{\otimes 2}} \simeq K_{\calX}$.

    The d\=/critical scheme $\Crit (f)$
    in~\cref{para-d-crit}~\cref{itm-d-crit-1}
    has a canonical orientation given by
    $K_{\Crit (f)}^{\smash{1/2}} = K_U |_{\Crit (f)^\red}$.

    By \cite[Lemma~2.58]{Joyce2015},
    for a smooth morphism $g \colon \calY \to \calX$ of d\=/critical stacks,
    compatible with the d\=/critical structures,
    an orientation~$K_{\calX}^{\smash{1/2}}$ of~$\calX$
    induces an orientation of~$\calY$ given by
    $K_{\calY}^{\smash{1/2}} = g^* (K_{\calX}^{\smash{1/2}}) \otimes \det \bbL_{\calY / \calX} |_{\calY^\red}$.
\end{para}

\begin{para}[Definition for schemes]
    Let $X$ be an oriented d\=/critical $\bbK$\=/scheme.
    Its \emph{motivic Behrend function}
    $\nu^\mot_X \in \Mm (X)$ is defined by the following property:

    \begin{itemize}
        \item
            For any critical chart $i \colon \Crit (f) \hookrightarrow X$,
            where $f \colon U \to \bbA^1$ and $U$ is a smooth $\bbK$\=/variety, we have
            \begin{equation}
                \label{eq-nu-chart}
                i^* (\nu^\mot_X) =
                -\bbL^{-{\dim U / 2}} \cdot
                \Phi_f ([U]) \cdot
                \Upsilon \bigl( i^* (K_X^{\smash{1/2}}) \otimes
                K_U^{-1} |_{\Crit (f)^\red} \bigr) \ ,
            \end{equation}
            in $\Mm (\Crit (f))$,
            where $\Phi_{\smash{f}}$ is the vanishing cycle map
            defined in~\cref{para-def-nearby},
            and $\Phi_{\smash{f}} ([U])$ is supported on~$\Crit (f)$.
            The map~$\Upsilon$ is as in~\cref{para-upsilon},
            and the part inside $\Upsilon ({\cdots})$ is a line bundle on~$\Crit (f)^\red$
            whose square is trivial, so it can be seen as a $\upmu_2$-bundle.
    \end{itemize}
    This is well-defined due to Bussi--Joyce--Meinhardt~%
    \cite[Theorem~5.10]{BJM2019}.

    For $X$ as above, and a smooth morphism $g \colon Y \to X$
    of relative dimension~$d$, where $Y$ is equipped with the induced
    oriented d\=/critical structure, we have the relation
    \begin{equation}
        \label{eq-nu-sm-pb-sch}
        g^* (\nu^\mot_X) =
        \bbL^{d / 2} \cdot
        \nu^\mot_Y \ ,
    \end{equation}
    which follows from \cite[Theorem~5.14]{BBBBJ2015}.
\end{para}

\begin{para}[Definition for stacks]
    \label{para-nu-mot-st}
    Let $\calX$ be an oriented d\=/critical stack over~$\bbK$,
    and assume that $\calX$ is Nisnevich locally a quotient stack
    in the sense of~\cref{para-quot-st}.

    We define the \emph{motivic Behrend function} of~$\calX$,
    slightly generalizing the construction of
    Ben-Bassat--Brav--Bussi--Joyce~%
    \cite[Theorem~5.14]{BBBBJ2015},
    who only considered stacks that are
    Zariski locally quotient stacks.
\end{para}

\begin{theorem*}
    Let $\calX$ be as above.
    Then there exists a unique element
    \begin{equation*}
        \nu^\mot_{\calX} \in \Mmhat (\calX) \ ,
    \end{equation*}
    called the \emph{motivic Behrend function} of $\calX$,
    such that for any $\bbK$\=/variety~$Y$
    and any smooth morphism $f \colon Y \to \calX$ of relative dimension~$d$,
    we have
    \begin{equation}
        \label{eq-nu-sm-pb-sch-stack}
        f^* (\nu^\mot_{\calX}) =
        \bbL^{d / 2} \cdot \nu^\mot_Y
    \end{equation}
    in $\Mmhat (Y)$, where $Y$ is equipped with the induced
    oriented d\=/critical structure.

    In particular, we write $\nu^\mot_{\frX} = \nu^\mot_{\calX}$
    if the d\=/critical structure on $\calX$ comes from
    a $(-1)$\=/shifted symplectic stack~$\frX$ with $\calX \simeq \frX_\mathrm{cl}$.
\end{theorem*}

\begin{proof}
    We first show that the theorem holds
    when $\calX = X$ is an algebraic space.
    Indeed, this follows formally from
    \cref{thm-nisn-desc} and
    the relation~\cref{eq-nu-sm-pb-sch} for schemes,
    since $X$ has a Nisnevich cover by affine varieties.

    Also, note that if the element $\nu^\mot_{\calX}$ exists,
    then the relation~\cref{eq-nu-sm-pb-sch-stack} must also hold
    for smooth morphisms from algebraic spaces~$Y$ to~$\calX$,
    by passing to a Nisnevich cover of~$Y$ by affine varieties.

    Now, the proof of \cite[Theorem~5.14]{BBBBJ2015}
    can be repeated word-by-word
    to show that the theorem is true when
    $\calX \simeq [S / G]$ is a quotient stack,
    where $S$ is an algebraic space over~$\bbK$
    and $G = \GL (n)$ for some~$n$.

    For the general case, let $(j_i \colon \calX_i \hookrightarrow \calX)_{i \in I}$
    be a Nisnevich cover by quotient stacks.
    The condition on $\nu^\mot_{\calX}$ forces that
    $j_i^* (\nu^\mot_{\calX}) = \nu^\mot_{\calX_i}$ for all~$i$.
    We show that the elements
    $\nu^\mot_{\calX_i}$ agree on overlaps.
    Indeed, let $1, 2 \in I$ be two indices,
    and let $\calX_{1, 2} = \calX_1 \times_{\calX} \calX_2$.
    Then $\calX_{1, 2}$ is also a quotient stack,
    so the theorem holds for~$\calX_{1, 2}$.
    Let $j'_i \colon \calX_{1, 2} \to \calX_i$ be the projections, where $i = 1, 2$.
    Then we have $(j'_i)^* (\nu^\mot_{\calX_i}) = \nu^\mot_{\calX_{1, 2}}$
    for $i = 1, 2$,
    since the left-hand side satisfies the
    characterizing property of~$\nu^\mot_{\calX_{1, 2}}$.
    By \cref{thm-nisn-desc},
    it then follows that the elements $\nu^\mot_{\calX_i}$ for~$i \in I$
    glue to a unique element~$\nu^\mot_{\calX}$,
    and a standard argument verifies that
    it satisfies the relation~\cref{eq-nu-sm-pb-sch-stack}.
\end{proof}

\begin{para}[Compatibility with smooth pullbacks]
    \label[theorem]{thm-nu-sm-pb}
    We now show that the smooth pullback relation~\cref{eq-nu-sm-pb-sch-stack}
    holds for all smooth morphisms of d\=/critical stacks.
\end{para}

\begin{theorem*}
    Let $\calX, \calY$ be oriented d\=/critical stacks over~$\bbK$
    that are Nisnevich locally quotient stacks,
    and let $f \colon \calY \to \calX$ be a smooth morphism
    of relative dimension~$d$ which is compatible with
    the oriented d\=/critical structures.
    Then we have the relation
    \begin{equation}
        \label{eq-nu-sm-pb}
        f^* (\nu^\mot_{\calX}) =
        \bbL^{d / 2} \cdot \nu^\mot_{\calY} \ .
    \end{equation}
\end{theorem*}

\begin{proof}
    It is straightforward to verify that the element
    $\bbL^{-d/2} \cdot f^* (\nu^\mot_{\calX})$
    satisfies the characterizing property of~$\nu^\mot_{\calY}$.
\end{proof}

\begin{para}[The numerical Behrend function]
    Let $\calX$ be an algebraic stack over~$\bbK$
    that is Nisnevich locally a quotient stack,
    equipped with an oriented d\=/critical structure.
    The \emph{Behrend function} of~$\calX$ is the constructible function
    \begin{equation*}
        \nu_{\calX} = \chi (\nu^\mot_{\calX}) \in \CF (\calX) \ ,
    \end{equation*}
    where $\chi$~denotes taking the pointwise Euler characteristic,
    as in~\cref{para-cf}.

    In fact, one can still define~$\nu_{\calX}$
    even if $\calX$ is only \'etale locally a quotient stack,
    and without the orientability assumption.
    Indeed, the relation~\cref{eq-nu-sm-pb}
    implies that the numerical Behrend function is compatible with
    smooth morphisms preserving the d\=/critical structure
    (not necessarily orientations),
    up to a sign~$(-1)^d$, where $d$~is the relative dimension.
    This is because changing the orientation only affects
    the term $\Upsilon ({\cdots})$ in~\cref{eq-nu-chart},
    which always has Euler characteristic~$1$,
    and the sign is due to the fact that $\chi (\bbL^{1/2}) = -1$.
    Now, to define~$\nu_{\calX}$, one can pass to a smooth cover of~$\calX$
    by $\bbK$\=/varieties, and apply smooth descent of constructible functions.

    When $\bbK = \bbC$, the Behrend function~$\nu_{\calX}$
    agrees with the original definitions by
    Behrend~\cite{Behrend2009} and Joyce--Song~\cite[\S4.1]{JoyceSong2012}.
    This follows from the compatibility
    of both versions with smooth pullbacks,
    namely \cref{thm-nu-sm-pb} and \cite[Theorem~4.3]{JoyceSong2012},
    and the case of critical loci on smooth varieties, which follows from
    \cite[Theorem~3.10]{DenefLoeser2002}
    and \cite[Theorem~4.7]{JoyceSong2012}.
\end{para}

%% file: grad.tex
We discuss the \emph{stack of graded points}
and the \emph{stack of filtered points}
of algebraic stacks and derived algebraic stacks,
following Halpern-Leistner~\cite{HalpernLeistnerInstability},
and study their interactions with shifted symplectic structures.
Then, in~\cref{subsec-grad-loc},
we give local descriptions of these stacks
using \'etale covers of the original stack by quotient stacks.

\subsection{Definition and deformation theory}
\label{subsec-grad-def}

\begin{para}[Definition]
    \label{para-def-grad}
    Let $\calX$ be an algebraic stack over~$\bbK$.
    Following Halpern-Leistner~\cite{HalpernLeistnerInstability},
    we define the \emph{stack of graded points}
    and the \emph{stack of filtered points} of $\frX$, respectively,
    as the mapping stacks over $\bbK$,
    \begin{align*}
        \Grad (\calX) & =
        \Map ([* / \Gm], \calX) \ ,
        \\
        \Filt (\calX) & =
        \Map ([\bbA^1 / \Gm], \calX) \ ,
    \end{align*}
    where $\Gm$ acts on $\bbA^1$ by scaling.

    Consider the morphisms
    \begin{equation*}
        \begin{tikzcd}
            {} [* / \Gm]
            \ar[shift left=0.5ex, r, "0"]
            \ar[loop, in=195, out=165, looseness=4, "(-)^{-1}"',
                start anchor={[yshift=1.25ex]west},
                end anchor={[yshift=-1.25ex]west}]
            &
            {} [\bbA^1 / \Gm]
            \ar[shift left=0.5ex, l, "\pr"]
            &
            * \vphantom{[* / \Gm]} \ ,
            \ar[shift left=0.5ex, l, "1"]
            \ar[shift right=0.5ex, l, "0"']
            \ar[ll, bend right, start anchor=north west, end anchor=north east, looseness=.8]
        \end{tikzcd}
    \end{equation*}
    where $\pr$ is induced by the projection $\bbA^1 \to *$.
    These morphisms induce morphisms of stacks
    \begin{equation*}
        \begin{tikzcd}
            \Grad (\calX)
            \ar[rr, bend left, start anchor=north east, end anchor=north west, looseness=.8, "\tot"]
            \ar[loop, in=195, out=165, looseness=4, "\mathrm{op}"',
                start anchor={[yshift=1.25ex]west},
                end anchor={[yshift=-1.25ex]west}]
            \ar[r, shift right=0.5ex, "\ssf"']
            &
            \Filt (\calX)
            \ar[l, shift right=0.5ex, "\gr"']
            \ar[r, shift left=0.5ex, "\ev_0"]
            \ar[r, shift right=0.5ex, "\ev_1"']
            &
            \calX \vphantom{\Grad (\calX)} \rlap{ ,}
        \end{tikzcd}
    \end{equation*}
    where the notations `$\op$', `$\gr$', `$\ssf$', and `$\tot$' stand for
    the \emph{opposite graded point},
    the \emph{associated graded point},
    the \emph{split filtration},
    and the \emph{total point}, respectively.

    The stacks $\Grad (\calX)$ and $\Filt (\calX)$ are algebraic stacks over~$\bbK$.
    If, moreover, $\calX$ has affine stabilizers,
    then so do $\Grad (\calX)$ and $\Filt (\calX)$,
    and by~\cite[Lemma~1.3.8]{HalpernLeistnerInstability},
    the morphism~$\gr$ is quasi-compact,
    and induces a bijection $\pi_0 (\Filt (\calX)) \simto \pi_0 (\Grad (\calX))$.
\end{para}

\begin{para}[The derived version]
    Now, consider a derived algebraic stack~$\frX$ over~$\bbK$.
    Similarly, we can consider the derived mapping stacks
    \begin{align*}
        \dGrad (\frX) & =
        \dMap ([* / \Gm], \frX) \ ,
        \\
        \dFilt (\frX) & =
        \dMap ([\bbA^1 / \Gm], \frX) \ ,
    \end{align*}
    where $\dMap (-, -)$ denotes the mapping stack
    in the $\infty$\=/category of derived stacks.
    These are again derived algebraic stacks over~$\bbK$,
    by Halpern-Leistner--Preygel~%
    \cite[Theorem~5.1.1]{HalpernLeistnerPreygel2023}
    or Halpern-Leistner~\cite[Theorem~1.2.1]{HalpernLeistnerDerived}.
    Moreover, they are locally finitely presented whenever $\frX$ is,
    by an argument similar to that in \cite[\S5.1.5]{HalpernLeistnerPreygel2023}.

    We use the same notations $\gr$, $\ev_1$, etc.,
    as in \cref{para-def-grad},
    for the induced morphisms between these derived stacks.

    We often write $\Grad (\frX), \Filt (\frX)$
    for $\dGrad (\frX), \dFilt (\frX)$
    when there is no ambiguity.
    We adopt the convention that $\Grad (\calX)$ and $\Filt (\calX)$
    refer to the classical versions when $\calX$ is assumed to be classical,
    and the derived versions otherwise.

    Note that for a classical algebraic stack~$\calX$,
    the stacks $\dGrad (\calX)$ and $\dFilt (\calX)$
    can have non-trivial derived structure.
    See~\cite[Example~1.6.4]{HalpernLeistnerDerived} for an example of this phenomenon.
    However, we always have $\dGrad (\calX)_\mathrm{cl} \simeq \Grad (\calX)$
    and $\dFilt (\calX)_\mathrm{cl} \simeq \Filt (\calX)$,
    where $(-)_\mathrm{cl}$ denotes the classical truncation.
    More generally, by~\cite[Lemma~1.2.1]{HalpernLeistnerInstability},
    for a derived algebraic stack~$\frX$,
    we have $\dGrad (\frX)_\mathrm{cl} \simeq \Grad (\frX_\mathrm{cl})$
    and $\dFilt (\frX)_\mathrm{cl} \simeq \Filt (\frX_\mathrm{cl})$.
\end{para}

\begin{para}[For quotient stacks]
    \label{para-grad-quot}
    Let $\calX = [U / G]$ be a quotient stack,
    where $U$ is an algebraic space over~$\bbK$,
    acted on by a reductive group~$G$.
    The stacks of graded and filtered points of~$\calX$
    can be described very explicitly,
    following Halpern-Leistner~\cite[\S1.4]{HalpernLeistnerInstability}.

    Let $\lambda \colon \Gm \to G$ be a \emph{cocharacter},
    that is, a morphism of algebraic groups.
    Define the \emph{Levi subgroup}
    and the \emph{parabolic subgroup} of~$G$
    associated to~$\lambda$ by
    \begin{align*}
        L_\lambda & =
        \{ g \in G \mid g = \lambda (t) \, g \, \lambda (t)^{-1} \text{ for all } t \} \ ,
        \\
        P_\lambda & =
        \{ g \in G \mid \lim_{t \to 0} \lambda (t) \, g \, \lambda (t)^{-1} \text{ exists} \} \ ,
    \end{align*}
    respectively. Define the
    \emph{fixed locus} and the \emph{attractor} associated to~$\lambda$ by
    \begin{align*}
        U^{\lambda, 0} & =
        \Map^{\Gm} (*, U) \ ,
        \\
        U^{\lambda, +} & =
        \Map^{\Gm} (\bbA^1, U) \ ,
    \end{align*}
    where $\Map^{\Gm} (-, -)$ denotes the $\Gm$-equivariant mapping space;
    $\Gm$ acts on $U$ via~$\lambda$, and on $\bbA^1$ by scaling.
    These are algebraic spaces by
    Drinfeld--Gaitsgory~\cite[Proposition~1.3.4 and Theorem~1.5.2]{DrinfeldGaitsgory2014}
    or Halpern-Leistner~\cite[Proposition~1.4.1]{HalpernLeistnerInstability}.
    There is a closed immersion $U^{\lambda, 0} \hookrightarrow U$,
    an unramified morphism $U^{\lambda, +} \to U$ given by evaluation at~$1$,
    and an affine morphism $U^{\lambda, +} \to U^{\lambda, 0}$ given by evaluation at~$0$.

    The $G$-action on~$U$ induce a $P_\lambda$-action on~$U^{\lambda, +}$
    and an $L_\lambda$-action on~$U^{\lambda, 0}$.
    Moreover, by \cite[Theorem~1.4.8]{HalpernLeistnerInstability},
    we have
    \begin{align*}
        \Grad (\calX) & \simeq
        \coprod_{\lambda \colon \Gm \to G} [U^{\lambda, 0} / L_\lambda] \ ,
        \\
        \Filt (\calX) & \simeq
        \coprod_{\lambda \colon \Gm \to G} [U^{\lambda, +} / P_\lambda] \ ,
    \end{align*}
    where the disjoint union is over all
    conjugacy classes of cocharacters~$\lambda$.
\end{para}

\begin{para}[Deformation theory]
    For a derived algebraic stack~$\frX$ locally of finite presentation over~$\bbK$,
    one can express the tangent complexes of $\Grad (\frX)$ and $\Filt (\frX)$
    in terms of that of~$\frX$. Concretely,
    by Halpern-Leistner--Preygel~\cite[Proposition~5.1.10]{HalpernLeistnerPreygel2023}
    or Halpern-Leistner~\cite[Lemma~1.2.2]{HalpernLeistnerInstability},
    we have
    \begin{align}
        \label{eq-grad-tan}
        \bbT_{\Grad (\frX)} & \simeq
        \tot^* (\bbT_{\frX})_0 \ ,
        \\
        \label{eq-filt-tan}
        \bbT_{\Filt (\frX)} & \simeq
        q_* \circ \ev^* (\bbT_{\frX}) \ ,
    \end{align}
    where $(-)_0$ denotes the weight~$0$ part
    with respect to the natural $\Gm$-action,
    $\ev \colon [\bbA^1 / \Gm] \times \Filt (\frX) \to \frX$
    is the evaluation morphism,
    and $q \colon [\bbA^1 / \Gm] \times \Filt (\frX) \to \Filt (\frX)$ is the projection.
\end{para}

\begin{para}[Oriented Lagrangian correspondences]
    Let $\frX$ and $\frY$ be oriented $s$\=/shifted symplectic stacks over~$\bbK$,
    as in \cref{para-ori}, where $s$~is odd, and let
    \begin{equation}
        \label{eq-def-lag-corr}
        \frX \overset{f}{\longleftarrow}
        \frL \overset{g}{\longrightarrow} \frY
    \end{equation}
    be an \emph{$s$\=/shifted Lagrangian correspondence},
    in the sense of \cite[\S2.4]{CalaqueHaugsengScheimbauer2022}.
    We thus have an exact triangle
    \begin{equation}
        \label{eq-lag-corr-tria}
        \bbT_{\frL} \longrightarrow
        f^* (\bbT_{\frX}) \oplus g^* (\bbT_{\frY}) \longrightarrow
        \bbL_{\frL} [s] \longrightarrow
        \bbT_{\frL} [1]
    \end{equation}
    of perfect complexes on~$\frL$.
    An \emph{orientation} of the
    shifted Lagrangian correspondence~\cref{eq-def-lag-corr}
    is an isomorphism
    $K_{\frL} \simeq f^* (K_{\frX}^{\smash{1/2}}) \otimes g^* (K_{\frY}^{\smash{1/2}})$,
    such that it squares to the canonical isomorphism
    $K_{\frL}^{\smash{\otimes 2}} \simeq f^* (K_{\frX}) \otimes g^* (K_{\frY})$
    induced by the exact triangle~\cref{eq-lag-corr-tria}.
\end{para}

\begin{theorem}
    \label{thm-filt-lag}
    Let $\frX$ be an $s$\=/shifted symplectic stack over~$\bbK$,
    with symplectic form~$\omega$.
    Then we have an induced $s$\=/shifted symplectic structure
    $\tot^* (\omega)$ on $\Grad (\frX)$,
    and an $s$\=/shifted Lagrangian correspondence
    \begin{equation}
        \label{eq-filt-lag}
        \Grad (\frX) \overset{\gr}{\longleftarrow}
        \Filt (\frX) \overset{\ev_1}{\longrightarrow} \frX \ .
    \end{equation}
    Moreover, if\/ $s$ is odd and $\frX$ has an orientation~%
    $K_{\frX}^{\smash{1/2}}$,
    then $\Grad (\frX)$ has an induced orientation~%
    $K_{\Grad (\frX)}^{\smash{1/2}}$,
    and the Lagrangian correspondence is oriented.
\end{theorem}

\begin{proof}
    The stacks $\Grad (\frX)$ and $\Filt (\frX)$ are derived algebraic stacks
    locally of finitely presentation over~$\bbK$,
    as mentioned in \cref{para-def-grad}.

    To prove that \cref{eq-filt-lag} is an
    $s$\=/shifted Lagrangian correspondence,
    by Calaque~\cite[Theorem~4.8]{Calaque2015},
    it is enough to show that the cospan
    \begin{equation}
        \label{eq-filt-or-cospan}
        [* / \Gm] \overset{0}{\longrightarrow}
        [\bbA^1 / \Gm] \overset{1}{\longleftarrow}
        *
    \end{equation}
    is a \emph{$0$-oriented cospan},
    in the sense of \cite[\S4.2]{Calaque2015}
    and \cite[\S2.5]{CalaqueHaugsengScheimbauer2022}.
    Indeed, $*$ carries a natural $0$-orientation,
    and the $0$-orientation on $[* / \Gm]$ is given by
    the isomorphism $\bbR \Gamma (\calO_{[* / \Gm]}) \simto \bbK$.
    To see that this is indeed a $0$-orientation,
    we check the condition in \cite[Definition~2.4]{PTVV2013}.
    For $A \in \cat{CdgA}_{\bbK}^{\smash{\leq 0}}$ and a perfect complex
    $\calE \in \Perf (\Spec A \times [* / \Gm])$,
    one has $p_* (\calE^\vee)^\vee \simeq p_* (\calE)$ on $\Spec A$,
    where $p \colon \Spec A \times [* / \Gm] \to \Spec A$ is the projection,
    since both sides are the weight~$0$ part of
    the induced $\Gm$-action on $\pi^* (\calE)$,
    where $\pi \colon \Spec A \to \Spec A \times [* / \Gm]$
    is the projection.

    To see that \cref{eq-filt-or-cospan} is a $0$-oriented cospan,
    we check the condition in
    \cite[Lemma~2.5.5]{CalaqueHaugsengScheimbauer2022}.
    For any $A \in \cat{CdgA}_{\bbK}^{\smash{\leq 0}}$
    and $\calE \in \Perf (\Spec A \times [\bbA^1 / \Gm])$,
    we need to show that the induced commutative diagram
    \begin{equation}
        \label{eq-or-cospan-diag}
        \begin{tikzcd}
            q_* (\calE) \ar[d] \ar[r] &
            p_* \circ 0^* (\calE) \ar[d] \\
            1^* (\calE) \ar[r] &
            q_* (\calE^\vee)^\vee
        \end{tikzcd}
    \end{equation}
    in $\Perf (A)$ is cartesian,
    where $p$ and $q$ are the projections from
    $\Spec A \times [* / \Gm]$ and
    $\Spec A \times [\bbA^1 / \Gm]$ to $\Spec A$, respectively.
    Indeed, as in
    Halpern-Leistner~\cite[Proposition~1.1.2~ff.]{HalpernLeistnerDerived},
    such an object $\calE$ can be seen as
    a filtered object in $\Perf (A)$, that is, a sequence of maps
    \begin{equation*}
        {\cdots} \longrightarrow
        E_{\geq 1} \longrightarrow
        E_{\geq 0} \longrightarrow
        E_{\geq -1} \longrightarrow
        {\cdots}
    \end{equation*}
    in $\Perf (A)$,
    where all but finitely many arrows are isomorphisms,
    such that $E_{\geq n} = 0$ for $n \gg 0$.
    Write $E_n = \operatorname{cofib} (E_{\geq n + 1} \to E_{\geq n})$,
    and write $E = \colim_{n \to -\infty} E_{\geq n}$.
    Then $0^* (\calE) \simeq \bigoplus_n E_n$, with the natural $\Gm$-action
    having weight $n$ on $E_n$.
    One can deduce from
    \cite[Proposition~1.1.2~ff.]{HalpernLeistnerDerived}
    that we have natural identifications
    \begin{align*}
        q_* (\calE) & \simeq E_{\geq 0} \ , \\
        p_* \circ 0^* (\calE) & \simeq E_0 \ , \\
        1^* (\calE) & \simeq E \ , \\
        q_* (\calE^\vee)^\vee & \simeq ((E^\vee)_{\geq 0})^\vee \simeq E_{\leq 0} \ ,
    \end{align*}
    where $E_{\leq 0} = \operatorname{cofib} (E_{\geq 1} \to E)$,
    and the arrows in the diagram~\cref{eq-or-cospan-diag}
    are the natural ones.
    This implies that~\cref{eq-or-cospan-diag} is cartesian.

    For the final statement,
    observe that
    \begin{align*}
        \tot^* (K_{\frX}) & \simeq
        \det \bigl( \tot^* (\bbL_{\frX})^0 \bigr) \otimes
        \det \bigl( \tot^* (\bbL_{\frX})^+ \bigr) \otimes
        \det \bigl( \tot^* (\bbL_{\frX})^- \bigr)
        \\ & \simeq
        K_{\Grad (\frX)} \otimes
        \det \bigl( \tot^* (\bbL_{\frX})^+ \bigr) \otimes
        \det \bigl( ( \tot^* (\bbL_{\frX})^+ )^\vee [-s] \bigr)
        \\ & \simeq
        K_{\Grad (\frX)} \otimes
        \det \bigl( \tot^* (\bbL_{\frX})^+ \bigr)^{\smash{2}} \ ,
    \end{align*}
    where $(-)^0$, $(-)^+$, $(-)^-$ denote the parts with
    zero, positive, and negative weights, respectively,
    with respect to the natural $\Gm$-action.
    Therefore, we may define
    \begin{equation}
        K_{\Grad (\frX)}^{\smash{1/2}} =
        \tot^* (K_{\frX}^{\smash{1/2}}) \otimes
        \det \bigl( \tot^* (\bbL_{\frX})^+ \bigr)^{-1} \ ,
    \end{equation}
    and this gives an orientation on $\Grad (\frX)$.
    To see that the $s$\=/shifted Lagrangian correspondence is oriented,
    consider the cartesian diagram
    \begin{equation}
        \label{eq-pf-lag-or}
        \begin{tikzcd}[column sep=small]
            \bbT_{\Filt (\frX)} \ar[r] \ar[d]
            \ar[dr, phantom, pos=.1, "\ulcorner"] &
            \gr^* (\bbT_{\Grad (\frX)}) \ar[d] \\
            \ev_1^* (\bbT_{\frX}) \ar[r] &
            \bbL_{\Filt (\frX)} [s]
        \end{tikzcd}
    \end{equation}
    in $\Perf (\Filt (\frX))$, witnessing the
    $s$\=/shifted Lagrangian correspondence structure.
    Write $\calE = \ev^* (\bbT_{\frX})$, where
    $\ev \colon [\bbA^1 / \Gm] \times \Filt (\frX) \to \frX$ is the evaluation morphism.
    As in the argument above,
    $\calE$ can be seen as a filtered object in $\Perf (\Filt (\frX))$,
    and the terms in~\cref{eq-pf-lag-or} can be identified with
    $E_{\geq 0}$, $E_0$, $E$, and $E_{\leq 0}$, respectively.
    In particular, one has
    $K_{\Filt (\frX)} \simeq \gr^* (K_{\Grad (\frX)}^{\smash{1/2}})
    \otimes \ev_1^* (K_{\frX}^{\smash{1/2}})$,
    as both sides can be identified with
    $\det (E_{\geq 0})^{-1}$.
\end{proof}

See also Kinjo--Park--Safronov~\cite[Corollary~3.18]{KinjoParkSafronov}
for a more detailed proof of the first part of this theorem.

\begin{lemma}
    \label{lem-filt-symp}
    Let $\frX$ be an $s$\=/shifted symplectic stack over~$\bbK$.
    Then we have an isomorphism
    \begin{equation*}
        \ssf^* (\bbT_{\Filt (\frX)}) \simeq
        \op^* \circ \ssf^* (\bbL_{\Filt (\frX)} [s])
    \end{equation*}
    of perfect complexes on~$\Grad (\frX)$.
\end{lemma}

\begin{proof}
    By Halpern-Leistner~\cite[Lemma~1.2.3]{HalpernLeistnerInstability},
    we have
    $\ssf^* (\bbT_{\Filt (\frX)}) \simeq \tot^* (\bbT_{\frX})_{\geq 0}$,
    where $(-)_{\geq 0}$ denotes taking the part with
    non-negative weights with respect to the natural $\Gm$-action.
    Consequently, we have
    $\op^* \circ \ssf^* (\bbT_{\Filt (\frX)}) \simeq \tot^* (\bbT_{\frX})_{\leq 0}$.
    Its dual shifted by~$s$ becomes
    $\tot^* (\bbL_{\frX} [s])_{\geq 0} \simeq \tot^* (\bbT_{\frX})_{\geq 0}$.
\end{proof}

\subsection{Local structure}
\label{subsec-grad-loc}

\begin{para}
    We now study the stacks of graded and filtered points
    of algebraic stacks
    admitting étale covers by quotient stacks
    in the sense of \crefrange{para-quot-st}{para-quot-st-eg},
    and give local descriptions of these stacks
    using such étale covers.
\end{para}

\begin{theorem}
    \label{thm-grad-filt-loc}
    Let $\calX$ be an algebraic stack over $\bbK$,
    and let $(\calX_i \to \calX)_{i \in I}$ be a representable étale cover,
    where each $\calX_i \simeq [S_i / G_i]$,
    with $S_i$ an algebraic space over~$\bbK$ and $G_i$ a reductive group.
    Then there are commutative diagrams
    \begin{equation}
        \label{eq-grad-filt-loc}
        \begin{tikzcd}[column sep={7em,between origins}]
            {} [S_i^{\lambda, 0} / L_{i, \lambda}] \ar[d] &
            {} [S_i^{\lambda, +} / P_{i, \lambda}] \ar[d]
            \ar[l] \ar[r] \ar[dl, phantom, pos=.2, "\urcorner"] &
            {} [S_i / G_i] \ar[d]
            \\
            \Grad (\calX) &
            \Filt (\calX) \ar[l, "\gr"'] \ar[r, "\ev_1"] &
            \calX \rlap{\textnormal{ ,}}
        \end{tikzcd}
    \end{equation}
    where all vertical arrows are representable and \'etale,
    $\lambda \colon \Gm \to G_i$ is a cocharacter,
    and the left-hand square is a pullback square.
    Moreover, the families
    \begin{align*}
        \bigl( [S_i^{\lambda, 0} / L_{i, \lambda}] & \longrightarrow \Grad (\calX) \bigr)
        _{i \in I, \ \lambda \colon \Gm \to G_i} \ ,
        \\
        \bigl( [S_i^{\lambda, +} / P_{i, \lambda}] & \longrightarrow \Filt (\calX) \bigr)
        _{i \in I, \ \lambda \colon \Gm \to G_i}
    \end{align*}
    are representable \'etale covers of
    $\Grad (\calX)$ and $\Filt (\calX)$, respectively.
\end{theorem}

\begin{proof}
    By Halpern-Leistner~\cite[Corollary~1.1.7]{HalpernLeistnerInstability},
    we have $\Grad (\calX_i) \simto \Grad (\calX) \times_{\calX} \calX_i$ for all~$i$.
    Therefore, the family $(\Grad (\calX_i) \to \Grad (\calX))_{i \in I}$
    is a representable étale cover of $\Grad (\calX)$.
    By \cref{lem-gr-et-pb} below,
    the family $(\Filt (\calX_i) \to \Filt (\calX))_{i \in I}$
    is a representable étale cover of $\Filt (\calX)$.
    The rest of the theorem follows from the description of
    $\Grad (\calX_i)$ and $\Filt (\calX_i)$ in \cref{para-grad-quot}.
    That the left-hand square in~\cref{eq-grad-filt-loc}
    is a pullback square follows from \cref{lem-gr-et-pb} below.
\end{proof}

\begin{para}[$\bbA^1$\=/action retracts]
    \label{para-a1-retract}
    We define a notion of \emph{$\bbA^1$\=/action retracts}
    for algebraic stacks,
    which are $\bbA^1$\=/deformation retracts
    that also gives rise to $\bbA^1$\=/actions.
    This will help us to prove \cref{lem-gr-et-pb} below.

    More precisely, for a morphism $\pi \colon \calY \to \calX$
    of algebraic stacks,
    a structure of a \emph{$\bbA^1$\=/action retract}
    consists of a monoid action $r \colon \bbA^1 \times \calY \to \calY$,
    a morphism $s \colon \calX \to \calY$,
    and an equivalence $r (0, -) \simeq s \circ \pi$.
    Here, the monoid structure on $\bbA^1$ is given by multiplication,
    and note that a monoid action requires extra coherence data.
    We say that $\pi$ is an \emph{$\bbA^1$\=/action retract}
    if such a structure exists.
\end{para}

\begin{lemma}
    \label{lem-gr-et-pb}
    Let $\calX, \calY$ be algebraic stacks over~$\bbK$ with affine stabilizers,
    and let $f \colon \calY \to \calX$ be an étale morphism.
    Then there is a pullback diagram
    \begin{equation}
        \label{eq-gr-et-pb}
        \begin{tikzcd}[column sep=1em]
            \Filt (\calY) \ar[r] \ar[d, "\gr"']
            \ar[dr, phantom, pos=.2, "\ulcorner"] &
            \Filt (\calX) \ar[d, "\gr"] \\
            \Grad (\calY) \ar[r] &
            \Grad (\calX) \rlap{\textnormal{ .}}
        \end{tikzcd}
    \end{equation}
\end{lemma}

\begin{proof}
    By Halpern-Leistner~\cite[Proposition~1.3.1]{HalpernLeistnerInstability},
    the horizontal arrows in~\cref{eq-gr-et-pb} are \'etale,
    and by \cite[Lemma~1.3.8]{HalpernLeistnerInstability},
    the vertical arrows are $\mathbb{A}^1$\=/action retracts.
    In particular, the induced morphism
    $\Filt (\calY) \to \Filt (\calX) \times_{\Grad (\calX)} \Grad (\calY)$
    is étale, and both sides are $\bbA^1$-action retracts
    onto $\Grad (\calY)$, so they are isomorphic by
    \cref{lem-et-gm} below.
\end{proof}

\begin{lemma}
    \label{lem-et-gm}
    Let $\calX, \calY_1, \calY_2$ be algebraic stacks over~$\bbK$, and let
    \begin{equation*}
        \begin{tikzcd}[column sep={3em,between origins}]
            \calY_1 \ar[rr, "f"] \ar[dr, "p_1"' pos=.4] &&
            \calY_2 \ar[dl, "p_2" pos=.4] \\
            & \calX
        \end{tikzcd}
    \end{equation*}
    be a commutative diagram,
    such that $p_1$ and $p_2$ are $\bbA^1$\=/action retracts,
    and~$f$ is $\mathbb{A}^1$-equivariant and étale.
    Then $f$ is an isomorphism.
\end{lemma}

\begin{proof}
    It is enough to show that $f$ is representable and bijective.
    Choose a point $y \in \calY_2 (\bbK)$, and form the pullback diagram
    \begin{equation*}
        \begin{tikzcd}
            \calU \ar[r, "g'"] \ar[d, "f'"']
            \ar[dr, phantom, pos=.2, "\ulcorner"] &
            \calY_1 \ar[d, "f"] \\
            \bbA^1 \ar[r, "g"] &
            \calY_2 \rlap{ ,}
        \end{tikzcd}
    \end{equation*}
    where~$g$ is given by $g (t) = t \cdot y$.
    Then~$\calU$ admits an $\bbA^1$-action retract to a point,
    and~$f'$ is $\bbA^1$-equivariant and étale.
    It is enough to show that~$f'$ is an isomorphism.

    The morphism~$f'$ is surjective since its image is a $\Gm$-equivariant
    open subset of~$\bbA^1$ containing $0$, and hence all of~$\bbA^1$.
    Let $u \in \calU (\bbK)$ be a point such that $f' (u) = 1$.
    Then the map $h \colon \bbA^1 \to \calU$ given by $h (t) = t \cdot u$
    is a section of~$f'$, and is hence also étale.
    Moreover, such a point~$u$ is unique,
    since another point~$u'$ with $f' (u') = 1$ would give a section
    $h' \colon \bbA^1 \to \calU$ with $h' (t) = t \cdot u'$,
    but $h$ and $h'$ agree on an open neighbourhood of~$0$, and hence everywhere.
    It follows that~$f'$ is bijective on~$\bbK$-points,
    and it is enough to show that~$\calU$ is an algebraic space.

    Consider the fibre product $V = \mathbb{A}^1 \times_{h, \calU, h} \mathbb{A}^1$,
    which is an algebraic space.
    Then~$V$ also admits an $\bbA^1$-action retract to a point,
    and projection $f'' \colon V \to \mathbb{A}^1$
    to either factor is étale and $\bbA^1$-equivariant.
    Repeating the above argument with~$V$ in place of~$\calU$,
    we see that $f''$ is bijective on~$\bbK$-points,
    and hence an isomorphism.
    It follows that $\bbK$-points of $\calU$
    cannot have non-trivial automorphisms.
\end{proof}

This lemma is based on a joint work~\cite{BuHalpernLeistnerIbanezNunezKinjo} with
D.~Halpern-Leistner, A.~Ibáñez Núñez, and T.~Kinjo,
and is an improvement of a similar result
in a previous version of this paper.

%% file: id.tex
We now present the main results of this paper,
in three different versions.
First, in~\cref{subsec-behrend-local},
we prove a local version of the main theorem,
using the theory of motivic nearby and vanishing cycles for stacks
developed in~\cref{subsec-nearby-stacks}.
Then, in~\cref{subsec-behrend-main},
we glue the local versions together
to prove the global version of our main result,
\cref{thm-behrend-main}.
Finally, in~\cref{subsec-behrend-num},
we take Euler characteristics in the main identity,
and obtain integral identities involving
the numerical Behrend functions.

\subsection{The local version}
\label{subsec-behrend-local}

\begin{theorem}
    \label{thm-behrend-local}
    Suppose that we are given the following data:

    \begin{itemize}
        \item
            A finite-dimensional $\Gm$-representation~$V$ over~$\bbK$. Let
            \begin{equation*}
                V = \bigoplus_{k \in \bbZ} V_k
            \end{equation*}
            be the decomposition into weight spaces.
            Write $V_+ = \bigoplus_{k > 0} V_k$.

        \item
            A $\bbK$\=/variety~$U$ acted on by~$\Gm$, and a
            $\Gm$-equivariant \'etale morphism $\iota \colon U \to V$.
            Write~$U^0 = U^{\Gm}$ for the fixed locus,
            and~$U^+ = \Map^{\Gm} (\bbA^1, U)$
            for the attractor as in \cref{para-grad-quot}.
            For a point $u_0 \in U^0 (\bbK)$, write
            \begin{equation*}
                U^+ (u_0) = \bigl\{ u \in U \bigm|
                \lim_{t \to 0} t \cdot u = u_0 \bigr\}
            \end{equation*}
            for the fibre of the limit map $U^+ \to U^0$ at~$u_0$.
            We have a canonical isomorphism $U^+ (u_0) \simeq V_+$
            by \cref{lem-et-gm}.

        \item
            A $\Gm$-invariant function $f \colon U \to \bbA^1$,
            with $f (u_0) = 0$.
    \end{itemize}
    Then we have the identities
    \begin{align}
        \label{eq-behrend-strong}
        \int \limits_{u \in U^+ (u_0)} \Psi_f ([U]) (u)
        & =
        \bbL^{\dim V_+} \cdot \Psi_{f} ([U^0]) (u_0) \ ,
        \\
        \label{eq-behrend-strong-vanishing}
        \int \limits_{u \in U^+ (u_0)} \Phi_f ([U]) (u)
        & =
        \bbL^{\dim V_+} \cdot \Phi_{f} ([U^0]) (u_0) \ .
    \end{align}
    Moreover, these hold as identities in~$\Mm (U^0)$,
    where we vary $u_0 \in U^0$.
\end{theorem}

This theorem can be seen as a generalization of
the integral identity conjectured by
Kontsevich--Soibelman~\cite[Conjecture~4]{KontsevichSoibelman2008},
and proved by Lê~\cite{Le2015},
who restricted to the case when the $\Gm$-action on $V$
only has weights $-1$, $0$, and~$1$.
Compare also Joyce--Song~\cite[Theorem~5.11]{JoyceSong2012},
where a similar identity involving Euler characteristics is proved,
with the same restriction on the weights.

The rest of this subsection is devoted to the proof of \cref{thm-behrend-local}.
In the following,
we first provide preliminaries on weighted projective spaces
and weighted blow-ups, and prove some preparatory results.
Then, in \cref{lem-behrend-weak},
we establish a weaker version of the theorem,
using the theory of motivic nearby cycles for stacks
developed in \cref{subsec-nearby-stacks}.
Finally, in \cref{para-pf-behrend-strong},
we show that the weaker version implies the stronger version.

\begin{para}[Weighted projective spaces]
    Let $V$ be a finite-dimensional $\Gm$-representation over~$\bbK$,
    with only positive weights.
    The \emph{weighted projective space} of $V$
    is the quotient stack
    \begin{equation*}
        \wP (V) = [(V \setminus \{ 0 \}) / \Gm] \ .
    \end{equation*}
    This is a proper Deligne--Mumford stack over~$\bbK$,
    since we have the identification
    \begin{equation*}
        \wP (V) \simeq \biggl[ \bbP (V) \bigg/
        \prod_{k=1}^{\dim V} \upmu_{n_k} \biggr] \ ,
    \end{equation*}
    where $\bbP (V)$ is the usual projective space,
    and using a basis of eigenvectors of~$V$,
    each $n_k$ is the weight of the $k$-th coordinate,
    and $\upmu_{n_k}$ acts by scaling the $k$-th coordinate.

    By \cref{para-pr-bun-rel}, the motive of~$\wP (V)$ is given by
    \begin{equation}
        [\wP (V)] = \frac{\bbL^{\dim V} - 1}{\bbL - 1} \ ,
    \end{equation}
    and is independent of the choice of weights on~$V$.

    We also consider the coarse space~$\cwP (V)$ of~$\wP (V)$,
    which is also given by
    \begin{equation*}
        \cwP (V) = \operatorname{Proj} \bbK [V] \ ,
    \end{equation*}
    where $\bbK [V]$ is the free polynomial algebra on~$V$,
    with $\bbZ$-grading given by the weights of~$V$.
    It is an integral, normal, projective $\bbK$\=/variety.
\end{para}

\begin{para}[Weighted blow-ups]
    \label{para-wbl}
    Let $V$ be a finite-dimensional $\Gm$-representation over~$\bbK$,
    with only positive weights.
    Let $U$ be a smooth $\bbK$\=/scheme, $U_0 \subset U$ a reduced closed subscheme,
    and let $p \colon U \to V$ be a smooth morphism
    such that $U_0 = p^{-1} (0)$.

    Define the \emph{weighted blow-up} of~$U$ along~$U_0$,
    with weights given by those of~$V$,
    as the quotient stack
    \begin{equation*}
        \wBl_{U_0} (U) =
        \Bigl[ \Bigl\{ (t, v, u) \in \bbA^1 \times (V \setminus \{ 0 \}) \times U \Bigm|
        p (u) = t \cdot v \Bigr\} \Big/ \Gm \Bigr] \ ,
    \end{equation*}
    where $t \cdot (-)$ denotes the $\Gm$-action
    naturally extended to~$t \in \bbA^1$,
    and $\Gm$ acts with weight~$-1$ on~$\bbA^1$,
    with the given weights on~$V$,
    and trivially on~$U$.
    Note that we have an isomorphism
    $\wBl_{U_0} (U) \simeq U \times_V \wBl_{\{ 0 \}} (V)$.

    The natural projection $\wBl_{U_0} (U) \to U$ is proper.
    It restricts to an isomorphism over $U \setminus U_0$,
    and has fibres $\wP (V)$ over points in~$U_0$.
    In particular, we have the relation
    \begin{equation}
        \label{eq-wbl-mot-rel}
        [\wBl_{U_0} (U)] = \frac{\bbL^{\dim V} - 1}{\bbL - 1} \cdot [U_0] +
        [U \setminus U_0]
    \end{equation}
    of motives on~$U$.
\end{para}

\begin{lemma}
    \label{lem-mot-torus-quot}
    Let $U$ be a separated algebraic space of finite type over~$\bbK$,
    acted on by a torus~$T \simeq \Gm^n$ for some~$n$,
    such that points in~$U$ have finite stabilizers.
    Let $\calX = [U / T]$ be the quotient stack.

    Then $\calX$ admits a coarse space $\pi \colon \calX \to X$
    which is a proper universal homeomorphism,
    and we have an isomorphism
    \begin{equation}
        \label{eq-mot-torus-quot}
        \pi_! = (\pi^*)^{-1} \colon \Mhat (\calX) \longsimto \Mhat (X) \ .
    \end{equation}
    A similar statement holds for $\Mm (\calX)$.
\end{lemma}

\begin{proof}
    Since~$U$ is separated and of finite type,
    the inertia $\calI_{\calX}$ is a closed substack
    of $H \times \calX$ for some finite group~$H \subset T$,
    and is thus finite over~$\calX$.
    It then follows from the Keel--Mori theorem~\cite{Conrad2005,KeelMori1997}
    that $\calX$ admits a coarse space $\pi \colon \calX \to X$,
    and that $\pi$ is a proper universal homeomorphism.

    To prove \cref{eq-mot-torus-quot},
    stratify $U$ by locally closed subspaces $U_i \subset U$,
    where each~$U_i$ is the locus of points
    with stabilizer~$H_i \subset T$,
    and let $\calX_i \subset \calX$ and $X_i \subset X$
    be the corresponding strata.
    Then $\calX_i \simeq (U_i / (T / H_i)) \times [* / H_i]$,
    so that $X_i \simeq U_i / (T / H_i)$ and
    $\calX_i \to X_i$ is a trivial gerbe with stabilizer~$H_i$.
    Note that the motive of $[* / \upmu_k]$ is~$1$ for all~$k > 0$,
    since it is the quotient $[\Gm / \Gm]$ with a weight~$k$ action,
    and $\Gm$ is a special group; see \cref{para-pr-bun-rel}.
    It follows that the motive of each~$[* / H_i]$ is~$1$,
    and the result follows.
\end{proof}

\begin{lemma}
    \label{lem-stab-pres}
    In the situation of \cref{thm-behrend-local},
    the locus in $U$ where the morphism~$\iota$
    preserves $\Gm$-stabilizers is open.
\end{lemma}

\begin{proof}
    For each $n > 1$, let $\zeta_n \in \Gm (\bbK)$ be a primitive $n$-th root of unity.
    It is enough to show that the locus of $u \in U$
    such that $\zeta_n \cdot u \neq u$ and
    $\iota (\zeta_n \cdot u) = \iota (u)$ is closed.
    The latter condition is equivalent to
    $\iota (u) \in V_\pn{n}$, where
    $V_\pn{n} = \bigoplus_{k \in \bbZ} V_{kn} \subset V$.
    Write $U_\pn{n} = \iota^{-1} (V_\pn{n})$,
    which is \'etale over $V_\pn{n}$,
    with a $\upmu_n$-action on its fibres,
    induced from the $\Gm$-action on~$U$.
    The locus where this action is trivial
    is open in~$U_\pn{n}$, proving the claim.
\end{proof}

\begin{lemma}
    \label{lem-git-quot-normal}
    In the situation of \cref{thm-behrend-local},
    suppose that $U$ is affine,
    and $\iota$ preserves $\Gm$-stabilizers
    and sends closed $\Gm$-orbits to closed $\Gm$-orbits.
    Then the affine GIT~quotient $U \slsl \Gm$ is normal.
\end{lemma}

\begin{proof}
    By Alper~\cite[Theorem~5.1]{Alper2013},
    since $\iota$ is \'etale and preserves $\Gm$-stabilizers,
    the induced morphism $\bar{\iota} \colon U \slsl \Gm \to V \slsl \Gm$
    is \'etale at $[u] \in (U \slsl \Gm) (\bbK)$
    for points $u \in U (\bbK)$ such that
    the $\Gm$-orbits of $u$ and $\iota (u)$ are closed.
    By the assumption on closed orbits, it is enough to require that
    the $\Gm$-orbit of $u$ is closed.
    Since every $S$-equivalence class in~$U$ contains a closed orbit,
    the morphism~$\bar{\iota}$ is \'etale, and it is enough to check that
    $V \slsl \Gm$ is normal.
    This follows from a standard fact in toric geometry,
    as in Cox--Little--Schenck \cite[Theorem~1.3.5]{CoxLittleSchenck2011},
    since $V \slsl \Gm \simeq \Spec \bbK [S]$
    for a saturated submonoid~$S \subset \bbZ^{\dim V}$.
\end{proof}

\begin{lemma}
    \label{lem-bij-iso}
    Let $f \colon X \to Y$ be a morphism of
    integral $\bbK$\=/varieties.
    If\/ $f$~is bijective on $\bbK$\=/points and $Y$~is normal,
    then $f$~is an isomorphism.
\end{lemma}

\begin{proof}
    By generic flatness and generic reducedness,
    $f$ is flat over a dense open subset $U \subset Y$
    with fibres $\Spec \bbK$, and hence \'etale,
    hence an isomorphism $f^{-1} (U) \simto U$.
    It follows that $f$ is birational.
    Now, a version of Zariski's main theorem
    \cite[IV{-}3, Corollary~8.12.10]{EGA}
    implies that $f$ is an open immersion, hence an isomorphism.
\end{proof}

\begin{lemma}
    \label{lem-behrend-weak}
    In the situation of \cref{thm-behrend-local},
    write $V_- = \bigoplus_{k < 0} V_k$,
    and for a point $u_0 \in U^0 (\bbK)$, consider the repeller
    \begin{equation*}
        U^- (u_0) = \bigl\{ u \in U \bigm|
        \lim_{t \to \infty} t \cdot u = u_0 \bigr\} \ ,
    \end{equation*}
    defined in the same way as $U^+ (u_0)$ for the opposite $\Gm$-action on~$U$.

    Then we have the identity
    \begin{equation}
        \label{eq-behrend-weak}
        \int \limits_{u \in U^+ (u_0)} \mspace{-9mu} \Psi_f ([U]) (u) -
        \int \limits_{u \in U^- (u_0)} \mspace{-9mu} \Psi_f ([U]) (u) =
        ( \bbL^{\dim V_+} - \bbL^{\dim V_-} ) \cdot
        \Psi_{f} ([U^0]) (u_0) \ .
    \end{equation}
    Moreover, this holds as an identity of monodromic motives on~$U^0$,
    where we vary $u_0 \in U^0$.
\end{lemma}

\begin{proof}
    Since $U$ is smooth,
    by Sumihiro~\cite[Corollary~2]{Sumihiro1974},
    $U$ admits a $\Gm$-invariant affine open cover.
    We may thus assume that $U$ is affine.
    Moreover, we apply this result whenever we shrink~$U$,
    so we may assume that $U$ is affine and connected throughout the proof.

    Write~$U^+, U^-$ for the attractor and repeller
    of the $\Gm$-action on~$U$.
    By Halpern-Leistner \cite[Propositions~1.3.1 and~1.3.2]{HalpernLeistnerInstability},
    the morphism $U^+ \to \iota^{-1} (V_+ \times V_0)$
    is étale and a closed immersion,
    and hence an open immersion.
    We may thus remove the closed subsets
    $\iota^{-1} (V_+ \times V_0) \setminus U^+$ and
    $\iota^{-1} (V_- \times V_0) \setminus U^-$ from~$U$,
    and assume that
    $U^\pm = \iota^{-1} (V_\pm \times V_0)$.
    The morphism~$\iota$ now sends closed $\Gm$-orbits to closed $\Gm$-orbits.

    By \cref{lem-stab-pres}, we may also assume that
    $\iota$ preserves $\Gm$-stabilizers,
    by replacing~$U$ with a $\Gm$-invariant open neighbourhood of~$U^0$.

    Let $U_\ominus = U \setminus U^-$, and let
    $U_\ominus^+ = U^+ \setminus U^0 \subset U_\ominus$.
    Consider the weighted blow-up
    \begin{equation*}
        \pi_\ominus \colon \tilde{U}_\ominus =
        \wBl_{U_\ominus^+} (U_\ominus) \longrightarrow U_\ominus \ ,
    \end{equation*}
    with weight $k$ along the $V_{-k}$-direction for $k > 0$,
    and write $\tilde{f}_\ominus = f \circ \pi_\ominus$.
    Explicitly, as in~\cref{para-wbl}, we may write
    \begin{align*}
        W_\ominus & = \Bigl\{ (t, v_-, u) \in \bbA^1 \times (V_- \setminus \{ 0 \}) \times U_\ominus \Bigm|
        \iota (u)_- = t^{-1} \cdot v_- \Bigr\} \ ,
        \\
        \tilde{U}_\ominus & = [W_\ominus / \Gm] \ ,
    \end{align*}
    where $\iota (u)_-$ is the projection of $\iota (u)$ to~$V_-$,
    and $\Gm$ acts on~$W_\ominus$ by
    $s \cdot (t, v_-, u) = (s^{-1} t, s^{-1} \cdot v_-, u)$.
    Note that $W_\ominus$ is smooth over~$\bbA^1 \times (V_- \setminus \{ 0 \})$,
    and hence over~$\bbK$.
    For any $u \in U_\ominus^+$, by
    \cref{thm-nearby-properties}~\cref{itm-proper-pushforward},
    we have
    \begin{align*}
        & \int \limits_{[v_-] \in \wP (V_-)} \Psi_{\tilde{f}_\ominus} ([\tilde{U}_\ominus]) ([v_-], u)
        \\[.5ex] = {} &
        \Psi_f ([\tilde{U}_\ominus]) (u)
        \\[.5ex] = {} &
        \Psi_f \bigl( [\wP (V_-) \times U_\ominus^+] + [U_\ominus \setminus U_\ominus^+] \bigr) (u)
        \\[.5ex] = {} &
        \bigl( [\wP (V_-)] - 1 \bigr) \cdot \Psi_f ([U_\ominus^+]) (u) +
        \Psi_f ([U_\ominus]) (u)
        \\ = {} &
        \biggl( \frac{\bbL^{\dim V_-} - 1}{\bbL - 1} - 1 \biggr) \cdot
        \Psi_f ([U_\ominus^+]) (u) +
        \Psi_f ([U]) (u) \ ,
        \numberthis
        \label{eq-pf-behrend-1}
    \end{align*}
    and this holds as an identity of monodromic motives on $U_\ominus^+$.

    Define $p^+ \colon U_\ominus^+ \to U^0$ by
    $p^+ (u) = \lim_{t \to 0} t \cdot u$.
    Then $f (u) = f (p^+ (u))$ for all $u \in U_\ominus^+$,
    and by \cref{thm-sch-nearby-properties}~\cref{itm-sch-smooth-pullback},
    we have
    \begin{equation}
        \label{eq-pf-behrend-2}
        \Psi_f ([U_\ominus^+]) (u) = \Psi_{f} ([U^0]) (p^+ (u))
    \end{equation}
    for all $u \in U_\ominus^+$.
    Again, this holds as an identity of monodromic motives on $U_\ominus^+$,
    where the right-hand side means $(p^+)^* \circ \Psi_f ([U^0])$.

    Now, consider the quotient stack
    \begin{equation}
        \check{U}_\ominus = [W_\ominus / \Gm^2] \ ,
    \end{equation}
    where $\Gm^2$ acts on~$W_\ominus$ by
    $(s_1, s_2) \cdot (t, v_-, u) =
    (s_1^{-1} t, s_1^{-1} s_2 \cdot v_-, s_2 \cdot u)$.
    There is, by definition, a principal $\Gm$-bundle
    $\tilde{\pi}_\ominus \colon \tilde{U}_\ominus \to \check{U}_\ominus$.
    There is a morphism $\check{f}_\ominus \colon \check{U}_\ominus \to \bbA^1$
    induced by~$\tilde{f}_\ominus$.

    Let $U \slsl \Gm$ be the affine GIT~quotient,
    and consider the reduced closed subscheme
    \begin{equation*}
        \tilde{U} \subset \cwP (V_+) \times \cwP (V_-) \times (U \slsl \Gm)
    \end{equation*}
    consisting of points $([\iota (u)_+], [\iota (u)_-], [u])$
    and $([v_+], [v_-], [u_0])$ for $u \in U$,
    $v_{\pm} \in V_{\pm} \setminus \{ 0 \}$, and $u_0 \in U^0$.
    There is a morphism $\tilde{f} \colon \tilde{U} \to \bbA^1$
    induced by~$f$.

    Consider the projection $\check{\pi}_\ominus \colon \check{U}_\ominus \to \tilde{U}$
    given by $(t, v_-, u) \mapsto ([\iota (u)_+], [v_-], [u])$.
    One can check that fibres of the composition $W_\ominus \to \tilde{U}$
    are single $\Gm^2$-orbits.
    We thus have an induced morphism
    $W_\ominus \slsl \Gm^2 \simto \tilde{U}$,
    which is an isomorphism by \cref{lem-bij-iso}.
    Here, we used the fact that $\tilde{U}$ is normal by \cref{lem-git-quot-normal},
    and the fact that $W_\ominus$ is integral since it is smooth and connected.
    In other words, the morphism $\check{\pi}_\ominus$
    is a coarse space map.
    In particular, it is proper by \cref{lem-mot-torus-quot}.

    Since the projection $\tilde{\pi}_\ominus \colon \tilde{U}_\ominus \to \check{U}_\ominus$
    is smooth and $\check{\pi}_\ominus$ is proper,
    by \cref{thm-nearby-properties} and
    \cref{lem-mot-torus-quot},
    for any $u \in U_\ominus^+$ and $[v_-] \in \wP (V_-)$, we have
    \begin{align*}
        \Psi_{\tilde{f}_\ominus} ([\tilde{U}_\ominus]) ([0, v_-, u])
        & =
        \Psi_{\check{f}_\ominus} ([\check{U}_\ominus]) ([0, v_-, u])
        \\ & =
        \Psi_{\tilde{f}} ([\tilde{U}]) ([\iota (u)_+], [v_-], [p^+ (u)]) \ ,
        \numberthis
        \label{eq-pf-behrend-3}
    \end{align*}
    where $[u] = [p^+ (u)]$ in $U \slsl \Gm$.
    Moreover, this holds as an identity of monodromic motives on
    $\wP (V_-) \times U_\ominus^+$.

    Combining \cref{eq-pf-behrend-1},
    \cref{eq-pf-behrend-2}, and \cref{eq-pf-behrend-3},
    we obtain the identity
    \begin{multline}
        \label{eq-pf-behrend-4}
        \Psi_f ([U]) (u) =
        \int \limits_{[v_-] \in \wP (V_-)} \Psi_{\tilde{f}} ([\tilde{U}]) ([\iota (u)_+], [v_-], [p^+ (u)])
        \\[-2ex]
        {} + \biggl( 1 - \frac{\bbL^{\dim V_-} - 1}{\bbL - 1} \biggr) \cdot
        \Psi_f ([U^0]) (p^+ (u)) \ ,
    \end{multline}
    where $u \in U_\oplus^+$ and $[v_-] \in \wP (V_-)$.
    Integrating over $u \in U^+ (u_0) \setminus \{ u_0 \}$, we obtain
    \begin{multline}
        \int \limits_{u \in U^+ (u_0) \setminus \{ u_0 \}} \Psi_f ([U]) (u) =
        (\bbL - 1) \cdot \int \limits_{([v_+], [v_-]) \in \wP (V_+) \times \wP (V_-)}
        \Psi_{\tilde{f}} ([\tilde{U}]) ([v_+], [v_-], [u_0])
        \\
        {} + ( \bbL^{\dim V_+} - 1 ) \cdot
        \biggl( 1 - \frac{\bbL^{\dim V_-} - 1}{\bbL - 1} \biggr) \cdot
        \Psi_f ([U^0]) (u_0) \ .
    \end{multline}
    Subtracting the analogous identity for
    integrating over $U_- (u_0) \setminus \{ u_0 \}$,
    we arrive at the desired identity~\cref{eq-behrend-weak}.
\end{proof}

\begin{para}[Proof of \cref{thm-behrend-local}]
    \label{para-pf-behrend-strong}
    Consider the $\Gm$-representation
    $V' = V \times \bbA^1$, with the $\Gm$-action on $V$ as given,
    and on $\smash{\bbA^1}$ by scaling.
    Let $U' = U \times \smash{\bbA^1}$,
    with the $\Gm$-action on $U$ as given,
    and on $\smash{\bbA^1}$ by scaling,
    and let $f' = f \circ \pr_1 \colon U' \to \bbA^1$,
    where $\pr_1 \colon U' \to U$ is the projection.
    Let $u'_0 = (u_0, 0) \in U'^0 = U^0 \times \{ 0 \}$.
    By \cref{thm-sch-nearby-properties}~\cref{itm-sch-smooth-pullback},
    we have $\Psi_{f'} ([U']) = \pr_1^* \circ \Psi_f ([U])$,
    and similarly, $\Psi_{f'} ([U'^0]) = \pr_1^* \circ \Psi_f ([U^0])$.

    Applying \cref{lem-behrend-weak} to this new set of data,
    and simplifying the expression by the observations above,
    we obtain
    \begin{equation*}
        \bbL \cdot \int \limits_{u \in U^+ (u_0)} \mspace{-9mu} \Psi_f ([U]) (u) -
        \int \limits_{u \in U^- (u_0)} \mspace{-9mu} \Psi_f ([U]) (u) =
        ( \bbL^{\dim V_+ + 1} - \bbL^{\dim V_-} ) \cdot
        \Psi_{f} ([U^0]) (u_0) \ .
    \end{equation*}
    Subtracting the original identity~\cref{eq-behrend-weak}
    from this, and dividing by $\bbL - 1$,
    we obtain the desired identity~\cref{eq-behrend-strong}.

    Finally, \cref{eq-behrend-strong-vanishing} follows from
    \cref{eq-behrend-strong} by the definition of $\Phi_f$.
    \qed
\end{para}

\subsection{The global version}
\label{subsec-behrend-main}

\begin{para}[Assumptions on the stack]
    \label{para-behrend-assumptions}
    In the following, we assume that
    $\frX$ is an oriented $(-1)$\=/shifted symplectic stack over~$\bbK$,
    with classical truncation $\calX = \frX_\mathrm{cl}$.

    We assume that $\calX$ is an algebraic stack that is
    Nisnevich locally fundamental
    in the sense of~\cref{para-quot-st}.
    For example, as in \cref{para-quot-st-eg},
    this is satisfied if $\calX$ admits a good moduli space,
    or can be covered by open substacks with good moduli spaces.
\end{para}

\begin{theorem}
    \label{thm-behrend-main}
    Let $\frX, \calX$ be as in~\cref{para-behrend-assumptions}.
    Consider the $(-1)$\=/shifted Lagrangian correspondence
    \begin{equation}
        \label{eq-filt-corr}
        \Grad (\frX) \overset{\gr}{\longleftarrow}
        \Filt (\frX) \overset{\ev_1}{\longrightarrow} \frX
    \end{equation}
    given by \cref{thm-filt-lag}.
    Then we have the identity
    \begin{equation}
        \label{eq-behrend-main}
        \gr_! \circ \ev_1^* (\nu^\mot_{\frX})
        = \bbL^{\vdim \Filt (\frX) / 2} \cdot
        \nu^\mot_{\Grad (\frX)}
    \end{equation}
    in $\Mmhat (\Grad (\calX))$,
    where $\vdim \Filt (\frX)$ is the virtual dimension of\/ $\Filt (\frX)$,
    seen as a function $\pi_0 (\Grad (\frX)) \simeq \pi_0 (\Filt (\frX)) \to \bbZ$.
\end{theorem}

We will prove the theorem in two steps.
First, in \cref{lem-behrend-nis},
we show that the theorem holds for a stack
if it holds for a Nisnevich cover of the stack,
reducing it to the case of fundamental stacks.
Then, we deduce the case of fundamental stacks
from the local version, \cref{thm-behrend-local}.

\begin{lemma}
    \label{lem-behrend-nis}
    Let $\frX, \calX$ be as in \cref{para-behrend-assumptions}.
    Let $(\calX_i \to \calX)_{i \in I}$
    be a Nisnevich cover,
    and write $\frX_i = \calX_i \times_{\calX} \frX$,
    with the induced $(-1)$\=/shifted symplectic structure
    and orientation.
    Then, if \cref{thm-behrend-main} holds for each~$\frX_i$,
    then it holds for~$\frX$.
\end{lemma}

\begin{proof}
    For each~$i$, consider the diagram
    \begin{equation}
        \begin{tikzcd}[column sep={6em,between origins}]
            \Grad (\calX_i) \ar[d] &
            \Filt (\calX_i) \ar[l, "\gr"'] \ar[r, "\ev_1"] \ar[d]
            \ar[dl, phantom, pos=.2, "\urcorner"] &
            \calX_i \ar[d] \\
            \Grad (\calX) &
            \Filt (\calX) \ar[l, "\gr"'] \ar[r, "\ev_1"] &
            \calX \rlap{ ,}
        \end{tikzcd}
    \end{equation}
    where the left-hand square is a pullback square
    by \cref{lem-gr-et-pb}.
    Therefore, there is a commutative diagram
    \begin{equation}
        \label{eq-cd-behrend-nis}
        \begin{tikzcd}[column sep={8em,between origins}]
            \Mmhat (\Grad (\calX_i)) &
            \Mmhat (\Filt (\calX_i)) \ar[l, "\gr_!"'] &
            \Mmhat (\calX_i) \ar[l, "\ev_1^*"'] \\
            \Mmhat (\Grad (\calX)) \ar[u] &
            \Mmhat (\Filt (\calX)) \ar[l, "\gr_!"'] \ar[u] &
            \Mmhat (\calX) \ar[l, "\ev_1^*"'] \ar[u] \rlap{ ,}
        \end{tikzcd}
    \end{equation}
    where the vertical maps are the pullback maps.

    By Halpern-Leistner~\cite[Corollary~1.1.7]{HalpernLeistnerInstability},
    we have $\Grad (\calX_i) \simto \Grad (\calX) \times_{\calX} \calX_i$ for all~$i$.
    Therefore, the family $(\Grad (\calX_i) \to \Grad (\calX))_{i \in I}$
    is a Nisnevich cover.
    By \cref{thm-nisn-desc}, it is enough to check
    the identity~\cref{eq-behrend-main}
    after pulling back to each~$\Grad (\calX_i)$.
    But this follows from the identity~\cref{eq-behrend-main} for each $\calX_i$,
    the commutativity of~\cref{eq-cd-behrend-nis},
    the relation~\cref{eq-nu-sm-pb} establishing
    the compatibility of the motivic Behrend function with smooth pullbacks,
    and the fact that the rank of the tangent complex of $\Filt (\frX_i)$
    agrees with that of $\Filt (\frX)$ on the corresponding components,
    which follows from~\cref{eq-filt-tan}.
\end{proof}

\begin{lemma}
    \label{lem-ori-compare}
    Suppose we have a pullback diagram of d\=/critical stacks
    \begin{equation}
        \begin{tikzcd}
            \calY' \ar[d, "f'"'] \ar[r, "g'"]
            \ar[dr, phantom, pos=.2, "\ulcorner"] &
            \calY \ar[d, "f"] \\
            \calX' \ar[r, "g"] &
            \calX \rlap{ ,}
        \end{tikzcd}
    \end{equation}
    where all morphisms are smooth
    and compatible with the d\=/critical structures.

    Let $K_{\calX}^{\smash{1/2}} \to \calX$ and $K_{\calY}^{\smash{1/2}} \to \calY$
    be orientations, not necessarily compatible with~$f$.
    Let $K_{\calX'}^{\smash{1/2}} \to \calX'$ and $K_{\calY'}^{\smash{1/2}} \to \calY'$
    be the orientations induced by
    $K_{\calX}^{\smash{1/2}}$ and $K_{\calY}^{\smash{1/2}}$,
    respectively, as in \cref{para-ori}.
    Then we have
    \begin{equation}
        g'^* \circ \Upsilon \bigl( K_{\calY}^{\smash{1/2}} \otimes
        f^* (K_{\calX}^{\smash{-1/2}}) \otimes
        \det (\bbL_{\calY / \calX})^{-1} \bigr)
        =
        \Upsilon \bigl( K_{\calY'}^{\smash{1/2}} \otimes
        f'^* (K_{\calX'}^{\smash{-1/2}}) \otimes
        \det (\bbL_{\calY' / \calX'})^{-1} \bigr)
    \end{equation}
    in $\Mmhat (\calY')$,
    where $\Upsilon$ is the map from \cref{para-upsilon},
    and the parts in $\Upsilon ({\cdots})$ are line bundles with trivial square,
    and can be seen as $\upmu_2$-bundles.
\end{lemma}

\begin{proof}
    These line bundles have trivial square
    by Joyce~\cite[Lemma~2.58]{Joyce2015}.
    We have
    \begin{align*}
        & \phantom{{} = {}}
        g'^* \bigl( K_{\calY}^{\smash{1/2}} \otimes
        f^* (K_{\calX}^{\smash{-1/2}}) \otimes
        \det (\bbL_{\calY / \calX})^{-1} \bigr)
        \\ & \simeq
        g'^* (K_{\calY}^{\smash{1/2}}) \otimes
        f'^* \circ g^* (K_{\calX}^{\smash{-1/2}}) \otimes
        \det (g'^* (\bbL_{\calY / \calX}))^{-1}
        \\ & \simeq
        K_{\calY'}^{\smash{1/2}} \otimes
        \det (\bbL_{\calY' / \calY})^{-1} \otimes
        f'^* (K_{\calX'}^{\smash{-1/2}}) \otimes
        f'^* \circ \det (\bbL_{\calX' / \calX}) \otimes
        \det (\bbL_{\calY' / \calX'})^{-1}
        \\ & \simeq
        K_{\calY'}^{\smash{1/2}} \otimes
        f'^* (K_{\calX'}^{\smash{-1/2}}) \otimes
        \det (\bbL_{\calY' / \calX'})^{-1} \rlap{ ,}
    \end{align*}
    and applying $\Upsilon$ gives the desired identity.
\end{proof}

\begin{para}[Proof of \cref{thm-behrend-main}]
    \label{para-pf-behrend-main}
    \allowdisplaybreaks
    By \cref{lem-behrend-nis}, we may assume that $\calX$ is fundamental.
    Let $\calX \simeq [S / G]$, where $S$ is an affine $\bbK$\=/variety,
    and $G = \GL (n)$ for some $n$.
    The classical truncation of
    the correspondence~\cref{eq-filt-corr} can be written as
    \begin{equation*}
        \coprod_{\lambda \colon \Gm \to G} {} [S^{\lambda, 0} / L_\lambda] \overset{\gr}{\longleftarrow}
        \coprod_{\lambda \colon \Gm \to G} {} [S^{\lambda, +} / P_\lambda] \overset{\ev_1}{\longrightarrow}
        [S / G] \ ,
    \end{equation*}
    with notations as in~\cref{para-grad-quot}.
    The assumption on~$G$ implies that all the groups
    $L_\lambda$ and $P_\lambda$ are special groups.

    We fix a cocharacter $\lambda \colon \Gm \to G$,
    and prove the identity on the component $[S^{\lambda, +} / P_\lambda]$.
    We may assume that $S^{\lambda, +} \neq \varnothing$.

    By Joyce~\cite[Remark~2.47]{Joyce2015},
    shrinking $S$ if necessary, we may assume that
    there exists a smooth affine $\bbK$\=/scheme $U$ acted on by $G$,
    and a $G$-invariant function $f \colon U \to \bbA^1$,
    such that $\calX$ is isomorphic as a d\=/critical stack
    to the critical locus $[\Crit (f) / G]$,
    and $S \simeq \Crit (f)$. We now have a commutative diagram
    \begin{equation}
        \begin{tikzcd}
            U^{\lambda, 0} \ar[d, "\pi^0"'] &
            U^{\lambda, +} \ar[l, "p"'] \ar[r, "i"] \ar[d, "\pi^+"] &
            U \ar[d, "\pi"] \\
            {} [U^{\lambda, 0} / L_\lambda] &
            {} [U^{\lambda, +} / P_\lambda] \ar[l, "\gr"'] \ar[r, "\ev_1"] &
            {} [U / G] \ .
        \end{tikzcd}
    \end{equation}

    Let $0 \in S^{\lambda, 0}$ be a $\bbK$\=/point,
    and let $V = \bbT_U |_0$ be the tangent space.
    Consider the $\Gm$-actions on $U$ and $V$ via the cocharacter~$\lambda$.
    By Luna~\cite[Lemma in \S III.1]{Luna1973},
    shrinking~$U$ if necessary,
    we may choose a $\Gm$-equivariant \'etale morphism
    $\iota \colon U \to V$ such that $\iota (0) = 0$.
    Applying \cref{thm-behrend-local} gives the identity
    \begin{equation}
        \label{eq-pf-behrend-main-loc}
        p_! \circ i^* \circ \Phi_f ([U]) =
        \bbL^{\dim V^\lambda_+} \cdot \Phi_f ([U^{\lambda, 0}]) \ ,
    \end{equation}
    where $V^\lambda_+ \subset V$ is the subspace
    where $\Gm$ acts with positive weights.
    Note that $\Phi_f (U)$ is supported on~$S$ by its definition.
    Let $K_S^{\smash{1/2}}$ be the orientation of
    the d\=/critical scheme~$S$ induced from that of~$\frX$.
    One computes that
    \begin{align*}
        \mathrlap{ \gr_! \circ \ev_1^* (\nu^\mot_{\frX}) }
        \quad &
        \\* & =
        [P_\lambda]^{-1} \cdot \gr_! \circ \pi^+_! \circ (\pi^+)^* \circ \ev_1^* (\nu^\mot_{\frX})
        \\ & =
        [P_\lambda]^{-1} \cdot \pi^0_! \circ p_! \circ i^* \circ \pi^* (\nu^\mot_{\frX})
        \\ & =
        \bbL^{\dim G / 2} \cdot [P_\lambda]^{-1} \cdot \pi^0_! \circ p_! \circ i^* (\nu^\mot_S)
        \\ & =
        -\bbL^{\dim G / 2 - \dim V / 2} \cdot
        [P_\lambda]^{-1} \cdot
        \pi^0_! \circ p_! \circ i^* \bigl(
            \Phi_f ([U]) \cdot \Upsilon (K_S^{\smash{1/2}} \otimes K_U^{-1} |_S)
        \bigr)
        \\ & =
        -\bbL^{\dim G / 2 - \dim V / 2} \cdot
        [P_\lambda]^{-1} \cdot {}
        \\* & \hspace{4em}
            \pi^0_! \circ p_! \Bigl(
                i^* \circ \Phi_f ([U]) \cdot
                i^* \circ \pi^* \circ \Upsilon \bigl( K_{\frX}^{\smash{1/2}} \otimes
                K_{[U / G]}^{-1} |_{\calX} \bigr)
            \Bigr)
        \\[.5ex] & =
        -\bbL^{\dim G / 2 - \dim V / 2} \cdot
        [P_\lambda]^{-1} \cdot {}
        \\* & \hspace{4em}
            \pi^0_! \circ p_! \Bigl(
                i^* \circ \Phi_f ([U]) \cdot
                (\pi^+)^* \circ \ev_1^* \circ
                \Upsilon \bigl( K_{\frX}^{\smash{1/2}} \otimes
                K_{[U / G]}^{-1} |_{\calX} \bigr)
            \Bigr)
        \\[.5ex] & =
        -\bbL^{\dim G / 2 - \dim V / 2} \cdot
        [P_\lambda]^{-1} \cdot {}
        \\* & \hspace{4em}
            \pi^0_! \circ p_! \Bigl(
                i^* \circ \Phi_f ([U]) \cdot
                (\pi^+)^* \circ \gr^* \circ
                \Upsilon \bigl( K_{\Grad (\frX)}^{\smash{1/2}} \otimes
                K_{[U^{\smash{\lambda, 0}} / L_\lambda]}^{-1} |_{[S^{\lambda, 0} / L_\lambda]} \bigr)
            \Bigr)
        \\[.5ex] & =
        -\bbL^{\dim G / 2 - \dim V / 2} \cdot
        [P_\lambda]^{-1} \cdot {}
        \\* & \hspace{4em}
            \pi^0_! \circ p_! \Bigl(
                i^* \circ \Phi_f ([U]) \cdot
                p^* \circ (\pi^0)^* \circ
                \Upsilon \bigl( K_{\Grad (\frX)}^{\smash{1/2}} \otimes
                K_{[U^{\smash{\lambda, 0}} / L_\lambda]}^{-1} |_{[S^{\lambda, 0} / L_\lambda]} \bigr)
            \Bigr)
        \\[.5ex] & =
        -\bbL^{\dim G / 2 - \dim V / 2} \cdot
        [P_\lambda]^{-1} \cdot {}
        \\* & \hspace{4em}
            \pi^0_! \Bigl(
                p_! \circ i^* \circ \Phi_f ([U]) \cdot
                (\pi^0)^* \circ
                \Upsilon \bigl( K_{\Grad (\frX)}^{\smash{1/2}} \otimes
                K_{[U^{\smash{\lambda, 0}} / L_\lambda]}^{-1} |_{[S^{\lambda, 0} / L_\lambda]} \bigr)
            \Bigr)
        \\[.5ex] & =
        -\smash{\bbL^{\dim G / 2 - \dim V / 2 + \dim V^\lambda_+}} \cdot
        [P_\lambda]^{-1} \cdot {}
        \\* & \hspace{4em}
            \pi^0_! \Bigl(
                \Phi_f ([U^{\lambda, 0}]) \cdot
                (\pi^0)^* \circ
                \Upsilon \bigl( K_{\Grad (\frX)}^{\smash{1/2}} \otimes
                K_{[U^{\smash{\lambda, 0}} / L_\lambda]}^{-1} |_{[S^{\lambda, 0} / L_\lambda]} \bigr)
            \Bigr)
        \\[.5ex] & =
        -\smash{\bbL^{\dim G / 2 - \dim V / 2 + \dim V^\lambda_+}} \cdot
        [P_\lambda]^{-1} \cdot {}
        \\* & \hspace{4em}
            \pi^0_! \Bigl(
                \Phi_f ([U^{\lambda, 0}]) \cdot
                \Upsilon \bigl( K_{S^{\lambda, 0}}^{\smash{1/2}} \otimes
                K_{U^{\lambda, 0}}^{-1} |_{S^{\lambda, 0}} \bigr)
            \Bigr)
        \\[.5ex] & =
        \smash{\bbL^{\dim G / 2 - \dim V / 2 + \dim V^\lambda_+ - \dim V^\lambda_0 / 2}} \cdot
        [P_\lambda]^{-1} \cdot
        \smash{\pi^0_! (\nu^\mot_{S^{\lambda, 0}})}
        \\ & =
        \smash{\bbL^{(\dim G - \dim L_\lambda) / 2 + (\dim V^\lambda_+ - \dim V^\lambda_-) / 2}} \cdot
        [P_\lambda]^{-1} \cdot
        \pi^0_! \circ (\pi^0)^* (\nu^\mot_{\Grad (\frX)})
        \\ & =
        \smash{\bbL^{(\dim G - \dim L_\lambda) / 2 + (\dim V^\lambda_+ - \dim V^\lambda_-) / 2}} \cdot
        [P_\lambda]^{-1} \cdot [L_\lambda] \cdot \nu^\mot_{\Grad (\frX)}
        \\* & =
        \smash{\bbL^{(\dim V^\lambda_+ - \dim V^\lambda_-) / 2}} \cdot
        \nu^\mot_{\Grad (\frX)} \ .
    \end{align*}
    Here, the first step uses~\cref{eq-pr-bun-rel};
    the third uses~\cref{eq-nu-sm-pb-sch-stack};
    the fourth uses~\cref{eq-nu-chart};
    the fifth uses \cref{lem-ori-compare},
    where the morphism~$f$ there is taken to be an isomorphism;
    the seventh uses the fact that
    the shifted Lagrangian correspondence~\cref{eq-filt-corr} is oriented,
    and the fact that the orientation for
    $\Grad ([\Crit (f) / G])$
    induced by the canonical one $K_{[U / G]}$
    is given by $K_{[U^{\lambda, 0} / L_\lambda]}$;
    the ninth uses~\cref{eq-mot-pf};
    the tenth is the key step, and uses~\cref{eq-pf-behrend-main-loc};
    the eleventh is analogous to the fifth;
    the twelfth uses~\cref{eq-nu-chart} again;
    the thirteenth uses~\cref{eq-nu-sm-pb-sch-stack} again;
    the fourteenth uses~\cref{eq-pr-bun-rel} again;
    and the final step uses the relation
    $[P_\lambda] = [L_\lambda] \cdot \bbL^{(\dim G - \dim L_\lambda) / 2}$.

    Finally, we verify that
    $\vdim \Filt^\lambda (\frX) = \dim V^\lambda_+ - \dim V^\lambda_-$,
    where $\Filt^\lambda (\frX) \subset \Filt (\frX)$ is the open and closed substack
    corresponding to the cocharacter~$\lambda$.
    Indeed, let $\frX' = [\Crit (f) / G]$ as a derived critical locus,
    with the natural $(-1)$\=/shifted symplectic structure,
    so $\frX'_{\smash{\mathrm{cl}}} \simeq \calX$.
    For $x \in S^{\lambda, 0} (\bbK)$,
    by \cref{lem-filt-symp}, one has
    \begin{align*}
        \rank (\bbL_{\Filt^\lambda (\frX)} |_x)
        & =
        \rank^{[0, 1]} (\bbL_{\Filt^\lambda (\frX)} |_x) -
        \rank^{[0, 1]} (\bbL_{\Filt^{-\lambda} (\frX)} |_x)
        \\* & =
        \rank^{[0, 1]} (\bbL_{\Filt^\lambda (\calX)} |_x) -
        \rank^{[0, 1]} (\bbL_{\Filt^{-\lambda} (\calX)} |_x)
        \\* & =
        \rank (\bbL_{\Filt^\lambda (\frX')} |_x) \ ,
        \numberthis
        \label{eq-vdim-filt-prime}
    \end{align*}
    where $\rank^{[0, 1]} = \dim \upH^0 - \dim \upH^1$.
    We have a presentation
    \begin{equation}
        \bbL_{\frX'} |_x \simeq \Bigl(
            \frg \longrightarrow
            \bbT_U |_{x} \longrightarrow
            \bbL_U |_{x} \longrightarrow
            \frg^\vee
        \Bigr)
    \end{equation}
    with degrees in $[-2, 1]$,
    where $\frg$~is the Lie algebra of~$G$.
    By Halpern-Leistner~\cite[Lemma~1.2.3]{HalpernLeistnerInstability},
    we have $\ssf^* (\bbL_{\smash{\Filt^\lambda (\frX)}}) \simeq \tot^* (\bbL_{\frX})_{\leq 0}$,
    where $(-)_{\leq 0}$ denotes the part of non-positive weights
    with respect to the natural $\Gm$-action. This now gives
    \begin{equation}
        \label{eq-tan-filt-quot}
        \bbL_{\Filt^\lambda (\frX')} |_x \simeq \Bigl(
            \frp_{\lambda} \longrightarrow
            \bbT_{U^{\lambda, -}} |_x \longrightarrow
            \bbL_{U^{\lambda, +}} |_x \longrightarrow
            \frp_{-\lambda}^\vee
        \Bigr) \ ,
    \end{equation}
    where $\frp_\lambda$ is the Lie algebra of~$P_\lambda$,
    and $-\lambda$ is the opposite cocharacter of~$\lambda$.
    Note that $\dim P_\lambda = \dim P_{-\lambda}$
    and that $\dim U^{\lambda, \pm} = \dim V^\lambda_{\pm} + \dim V^\lambda_0$.
    It follows that $\vdim \Filt (\frX)$,
    which is equal to the rank of~\cref{eq-tan-filt-quot}
    by~\cref{eq-vdim-filt-prime},
    is $\dim V^\lambda_+ - \dim V^\lambda_-$.
    \qed
\end{para}

%

\subsection{The numerical version}
\label{subsec-behrend-num}

\begin{para}[Assumptions on the stack]
    \label{para-benrend-num-assumptions}
    In the following, we assume that $\frX$ is a
    $(-1)$\=/shifted symplectic stack over~$\bbK$,
    with classical truncation $\calX = \frX_\mathrm{cl}$.
    Note that we no longer assume that $\frX$ is oriented.

    We assume that $\calX$ is an algebraic stack that is
    \'etale locally fundamental in the sense of~\cref{para-quot-st}.
    For example, as mentioned in~\cref{para-quot-st-eg},
    this is satisfied if $\calX$ has affine stabilizers
    and has reductive stabilizers at closed points.
\end{para}

\begin{para}
    For a graded point $\gamma \in \Grad (\calX) (\bbK)$, write
    \begin{equation*}
        \bbP (\gr^{-1} (\gamma)) = \Bigl(
            [* / \Gm] \underset{\Grad (\calX)}{\times} \Filt (\calX)
        \Bigr) \mathbin{\Big\backslash} \{ \ssf (\gamma) \} \ ,
    \end{equation*}
    where the map $[* / \Gm] \to \Grad (\calX)$ is given by
    the tautological $\Gm$-action on $\gamma$.
    The $\bbK$\=/point $\ssf (\gamma)$ is closed in the fibre product,
    which can be seen from
    the \'etale local description in~\cref{thm-grad-filt-loc},
    and $\{ \ssf (\gamma) \}$ denotes the corresponding closed substack.
    The space $\bbP (\gr^{-1} (\gamma))$ can be seen as
    the projectivized space of filtrations of a given
    associated graded point.

    As a remark, as mentioned in the proof of \cref{lem-gr-et-pb},
    the morphism $\gr$ is an $\bbA^1$\=/action retract, so the fibre
    $\gr^{-1} (\gamma) = \Spec (\bbK) \times_{\Grad (\calX)} \Filt (\calX)$
    is an $\bbA^1$-action retract to the point $\ssf (\gamma)$.
    The stack $\bbP (\gr^{-1} (\gamma))$ is the quotient
    of $\gr^{-1} (\gamma) \setminus \{ \ssf (\gamma) \}$
    by the $\Gm$-action which is part of this $\bbA^1$-action.
\end{para}

\begin{theorem}
    \label{thm-behrend-num}
    Let $\frX, \calX$ be as in \cref{para-benrend-num-assumptions}.
    Let $\gamma \in \Grad (\calX) (\bbK)$ be a graded point,
    and let $\bar{\gamma} = \op (\gamma)$ be its opposite graded point.

    Then we have the numerical identities
    \begin{align*}
        \nu_{\calX} (\tot (\gamma))
        & =
        (-1)^{\rank^{[0, 1]} ( \bbL_{\Filt (\calX)} |_{\ssf(\gamma)} ) -
            \rank^{[0, 1]} ( \bbL_{\Filt (\calX)} |_{\ssf(\bar{\gamma})} )} \cdot
        \nu_{\Grad (\calX)} (\gamma) \ ,
        \numberthis
        \label{eq-behrend-num-1}
        \\[1ex]
        \hspace{6em}
        & \hspace{-6em}
        \int \limits_{\varphi \in \bbP (\gr^{-1} (\gamma))} \nu_{\calX} (\ev_1 (\varphi)) \, d \chi -
        \int \limits_{\varphi \in \bbP (\gr^{-1} (\bar{\gamma}))} \nu_{\calX} (\ev_1 (\varphi)) \, d \chi
        \\
        & =
        \bigl( \dim \upH^0 (\bbL_{\Filt (\calX)} |_{\ssf(\gamma)})
        - \dim \upH^0 (\bbL_{\Filt (\calX)} |_{\ssf(\bar{\gamma})}) \bigr) \cdot
        \nu_{\calX} (\tot (\gamma)) \ ,
        \numberthis
        \label{eq-behrend-num-2}
    \end{align*}
    where $\rank^{[0, 1]} = \dim \upH^0 - \dim \upH^1$.
\end{theorem}

This theorem is a generalization of
Joyce--Song~\cite[Theorem~5.11]{JoyceSong2012},
who considered the case when $\calX$ is the moduli stack
of objects in a $3$-Calabi--Yau abelian category.

\begin{para}[Proof of \cref{thm-behrend-num}]
    \label{para-pf-behrend-num}
    By a similar argument as in the proof of \cref{lem-behrend-nis},
    passing to a representable \'etale cover of~$\calX$ by fundamental stacks,
    which induces representable \'etale covers of
    $\Grad (\calX)$ and $\Filt (\calX)$ by \cref{thm-grad-filt-loc},
    it is enough to prove the theorem when $\calX \simeq [S / G]$ is fundamental,
    where $S$ is an affine $\bbK$\=/scheme acted on by a reductive group~$G$.
    Here, we are using \'etale descent for constructible functions,
    instead of Nisnevich descent for rings of motives.

    As in \cref{para-pf-behrend-main},
    shrinking $S$ if necessary, we may assume that
    there exists a smooth affine $\bbK$\=/scheme $U$ acted on by $G$,
    and a $G$-invariant function $f \colon U \to \bbA^1$,
    such that $\calX$ is isomorphic as a d\=/critical stack
    to the critical locus $[\Crit (f) / G]$.
    Now, $\calX$ comes with a natural orientation,
    and the motivic Behrend function $\nu^\mot_{\calX}$ is defined.

    Applying \cref{thm-behrend-main},
    then evaluating the Euler characteristics at~$\gamma$,
    we obtain the identity
    \begin{equation}
        \label{eq-pf-behrend-num-1}
        \int \limits_{\varphi \in \gr^{-1} (\gamma)} \nu_{\calX} (\ev_1 (\varphi)) \, d \chi =
        (-1)^{\vdim_{\gamma} \Filt (\frX)} \cdot \nu_{\Grad (\calX)} (\gamma) \ .
    \end{equation}
    Let $\varphi_0 = \ssf (\gamma)$.
    Then the left-hand side of~\cref{eq-pf-behrend-num-1}
    is equal to $\nu_{\calX} (\ev_1 (\varphi_0)) = \nu_{\calX} (\tot (\gamma))$,
    since the integrand is $\Gm$-invariant
    and $\varphi_0$ is in the closure of all $\Gm$-orbits.
    Also, by \cref{lem-filt-symp}, we have
    \begin{equation}
        \vdim_\gamma \Filt (\frX) =
        \rank^{[0, 1]} \bbL_{\Filt (\calX)} |_{\ssf(\gamma)} -
        \rank^{[0, 1]} \bbL_{\Filt (\calX)} |_{\ssf(\bar{\gamma})} \ .
    \end{equation}
    This verifies~\cref{eq-behrend-num-1}.

    For~\cref{eq-behrend-num-2},
    apply \cref{thm-behrend-main} again,
    then take the difference of the evaluations
    at $\gamma$ and~$\bar{\gamma}$.
    This gives the identity
    \begin{multline}
        \label{eq-pf-behrend-num-2}
        (\bbL - 1) \cdot \biggl[
            \int \limits_{\varphi \in \bbP (\gr^{-1} (\gamma))} \nu^\mot_{\frX} (\ev_1 (\varphi)) -
            \int \limits_{\varphi \in \bbP (\gr^{-1} (\bar{\gamma}))} \nu^\mot_{\frX} (\ev_1 (\varphi))
        \biggr] \\[.5ex]
        {} +
        \bbL^{\dim \upH^1 (\bbL_{\Grad (\calX)} |_{\gamma})} \cdot
        \Bigl(
            \bbL^{-{\dim \upH^1 (\bbL_{\Filt (\calX)} |_{\ssf(\gamma)})}}
            - \bbL^{-{\dim \upH^1 (\bbL_{\Filt (\calX)} |_{\ssf(\bar{\gamma})})}}
        \Bigr) \cdot
        \nu^\mot_{\frX} (\tot (\gamma)) \\[1ex]
        = \Bigl(
            \bbL^{\rank (\bbL_{\Filt (\frX)} |_{\ssf(\gamma)}) / 2}
            - \bbL^{-{\rank (\bbL_{\Filt (\frX)} |_{\ssf(\gamma)})} / 2}
        \Bigr) \cdot
        \nu^\mot_{\Grad (\frX)} (\gamma)
    \end{multline}
    of monodromic motives over~$\bbK$.
    Here, we used the fact that the stabilizer group~$G_{\smash{\gamma}}$ of~$\gamma$
    in $\smash{\gr^{-1} (\gamma)}$ is special
    and has motive $\bbL^{\smash{\dim G_\gamma}}$,
    since $G_{\smash{\gamma}}$ is a subgroup of
    the fibre of the projection $P_\lambda \to L_\lambda$,
    and can be obtained by repeated extensions of $\Ga$.
    All of this can be seen by, for example,
    equivariantly embedding $S$ into an affine space with a linear $G$-action.

    Starting from~\cref{eq-pf-behrend-num-2},
    we divide both sides by $\bbL - 1$,
    and then take the Euler characteristic,
    which sets~$\bbL^{1/2}$ to~$-1$.
    We then apply the identity~\cref{eq-behrend-num-1}
    to convert $\nu_{\Grad (\calX)} (\gamma)$ to $\nu_{\calX} (\tot (\gamma))$.
    This gives the desired identity~\cref{eq-behrend-num-2}.
    \qed
\end{para}